\documentclass[journal=jacsat,manuscript=article]{achemso}

\usepackage{chemformula} % Formula subscripts using \ch{}
\usepackage[T1]{fontenc} % Use modern font encodings
\usepackage{amssymb} % for real numbers symbol
\usepackage{todonotes}
\usepackage{amsmath}
\usepackage{graphicx}
\usepackage{subcaption}
\usepackage{makecell}
\usepackage[final]{changes}
\usepackage{tikz}
\usepackage{tikz-3dplot}
\tdplotsetmaincoords{50}{125} % viewing angle (adjust as desired)

\newcommand{\be}{\begin{eqnarray}}
\newcommand{\ee}{\end{eqnarray}}

\def\refeq#1{(\ref{#1})}

\def\l{\left}
\def\r{\right}
\definechangesauthor[name={Martin Wlotzka}, color=blue]{MW}
\definechangesauthor[name={Norbert Asprion}, color=red]{NA}
\definechangesauthor[name={Pascal Schaefer}, color=purple]{PasS}
\definechangesauthor[name={Reviewer 1}, color=blue]{R1}
\definechangesauthor[name={Reviewer 2}, color=purple]{R2}
\definechangesauthor[name={Reviewers 1&2}, color=red]{R12}

\author{Jan Schwientek}
\author{Katrin Teichert}
\author{Jan Schröder}
\author{Johannes Höller}
\author{Patrick Schwartz}
\affiliation
{Fraunhofer ITWM, Fraunhofer-Platz 1, 67663 Kaiserslautern, Germany.}

\author{Norbert Asprion}
\author{Pascal Schäfer}
\author{Martin Wlotzka}
\affiliation
{BASF SE, Carl-Bosch-Strasse 38, 67056 Ludwigshafen am Rhein, Germany}

\author{Michael Bortz}
\affiliation
{Fraunhofer ITWM, Fraunhofer-Platz 1, 67663 Kaiserslautern, Germany.}

\email{Jan.Schwientek@itwm.fraunhofer.de}
\title
  {On Computing and Pricing of\linebreak Adjustable Robust Chemical Process Designs}


\begin{document}

\begin{abstract}
Model-based process simulation can be used to derive designs and operating conditions of chemical processes that optimally balance multiple objectives, such as quality, costs, or environmental impacts. This work focuses on identifying designs that hedge against uncertainties in model parameters to ensure feasibility, taking the possibility to adjust operating conditions into account. An adaptive scheme is proposed to pinpoint the relevant scenarios in a discretized uncertainty space; these scenarios are then fed into a multi-objective adjustable robust optimization framework reducing the computational burden compared to the consideration of all potential scenarios. Furthermore, we propose a method to quantify the cost or price of robustness, i.e., the compromise which has to be made in comparison to the nominal design case in order to hedge against uncertainty. The conceptual findings are illustrated with an industrially relevant case study.
\end{abstract}

\section{Introduction}

% Marking changes
%This is \added[id=auth,remark={Important addition}]{newly added text}
%This is \deleted[id=auth,remark={No longer needed}]{obsolete text}
%This is \replaced[id=auth]{updated}{outdated} content

Modeling of chemical processes is indispensable to identify and compare different plant designs and operating conditions according to key performance indicators \cite{bortzBook}. One such key performance indicator is robustness, i.e., the guarantee that despite uncertainties, e.g., in operating conditions or physical properties, relevant constraints are still fulfilled. These constraints may for example refer to product quality or safety requirements. For brevity, the term \emph{process design} will be used in the following to refer to both, plant designs and operating conditions.

\replaced[id=R1]{
For a long time now, the often possible distinction between \emph{here-and-now variables (HNV} or 1st-stage decisions), which have to be fixed before any uncertainty is observed, and \emph{wait-and-see variables (WSV) or 2nd-stage decisions}, which can be fixed after observing the realization of uncertainty, has been exploited to avoid overly pessimistic decisions. Initially, this approach was elaborated in two-stage stochastic optimization\cite{Birge2011} and later transferred to robust optimization\cite{BenTal2004}. The benefits in process engineering were demonstrated, for example, in \citeauthor{Lappas2016}\cite{Lappas2016} and also in a recent work of the authors\cite{schwientek2025}.
}{In a recent work \cite{schwientek2025}, it is shown that in order to avoid overly pessimistic process conditions, it is beneficial to distinguish between here-and-now and wait-and-see variables.} In process system engineering, the HNV refer to the plant design and include, for example, the kind and size of unit operations, while the WSV contain operational degrees of freedom like the reflux ratio in a distillation column or the temperature in a chemical reactor. As opposed to HNV, WSV are adjustable during process operation, i.e., they can be varied to react on uncertainties. Taking this adjustability into account allows to identify process conditions that guarantee the fulfillment of constraints while still offering the possibility to react to uncertainties.

\replaced[id=R1]{
Uncertainties are often of continuous nature. However, for computability reasons, these can only be considered in this way for certain models, e.g., where the uncertainties enter linearly in the objectives and constraints. Under general conditions as we have, a finite set of reference scenarios is often used to be able to calculate distributions in the stochastic or worst-case solutions in the robust approach, which is our concern, at all. Then, of course, those solutions are only approximations to the real ones.
}{Typically, a fixed set of discrete scenarios is chosen a priori to represent the continuous uncertainty set \cite{schwientek2025}.} Because the numerical effort scales with the number of scenarios, applying this approach to large-scale processes becomes computationally expensive or even intractable. This is addressed in this paper by introducing an adaptive \emph{worst-case\deleted[id=R2]{(WC)}} scenario selection built on a reference discretization. The adaptive scheme reduces the scenario count; the discretization helps to avoid globalization issues during optimization. Building on both, we propose a novel method for estimating (and visualizing) the price of robustness in a \emph{multi-criteria optimization (MO)} setting that accounts for the adjustability of the WSV.
%The same framework enables inverse robustness: the range of scenarios covered by the operating conditions is added as an extra objective in the MO, yielding precise trade-offs between robustness and key performance indicators such as cost efficiency, yield, and product quality.

The necessity of addressing uncertainties in the design of chemical processes has already been highlighted decades ago in \citeauthor{grossmann1978}\cite{grossmann1978}. Their work adopts a stochastic programming approach by weighting the scenarios. A more comprehensive and mathematically in-depth overview of the use of stochastic optimization techniques for the operation of chemical processes under uncertainty can be found in \citeauthor{henrion2001}\cite{henrion2001}. Consideration of uncertainties in model-based chemical process design also entails the creation of specialized software modules \cite{steimel2015}.

A significant focus within the robust approach has been the flexibility index, which quantifies the potential variability across different scenarios \cite{swaney1985part1,swaney1985part2}. \added[id=R1]{A discussion on the strong connection between flexibility analysis and two-stage robust optimization can be found, e.g., for the linear case, in \citeauthor{Zhang2016}\cite{Zhang2016}. Furthermore, the recent work \citeauthor{Zhang2026}\cite{Zhang2026} proposes an approach based on adjustable robust optimization with decision-dependent uncertainty, in which flexibility is directly integrated into the objective function via several measures in order to analyze the trade-off between design costs and operational flexibility when designing industrial systems.}

Regarding multi-criteria optimization under uncertainties, the work of \citeauthor{plazoglu1987}\cite{plazoglu1987}  was pioneering in integrating robustness aspects into an MO framework. There, a two-stage methodology is proposed that decomposes the MO problem into a series of structurally fixed single-objective problems using scalarization techniques. Similarly, in \citeauthor{dantus1999}\cite{dantus1999} the trade-off between costs and environmental effects is explored through a stochastic optimization method (simulated annealing), utilizing Monte-Carlo simulations to factor uncertainties across various process evaluations.

A more generalized version of the two-stage methodology was later introduced by \citeauthor{fu2000}\cite{fu2000}, where a master MO algorithm outlines a sequence of single-objective problems that incorporate uncertainty via sampling techniques, such as Hammersley sequence sampling. Subsequently, a comparable framework, including specialized genetic algorithms, was developed in \citeauthor{kheawhom2004,kheawhom2005}\cite{kheawhom2004,kheawhom2005}. This approach was further expanded by \citeauthor{kim2002}\cite{kim2002} to address solvent design challenges, where discrete optimization variables specify solvent structures. That research underscored how the probability distribution of uncertain parameters influences optimization results.

In addition to the two-stage framework, 
\citeauthor{chakraborty2003part2,chakraborty2003part3}\cite{chakraborty2003part2,chakraborty2003part3} implemented MO under uncertainty for plant-wide optimization using superstructures. They introduced a flexibility index to evaluate design robustness, applying it to analyze the Pareto boundary derived from an MO problem that balances costs against global warming potential. In a similar vein, \citeauthor{hoffmann2004}\cite{hoffmann2004} first computed a Pareto set based on deterministic assumptions, then examined the effects of parameter uncertainties by sampling their distributions and visualizing variations in objectives.

Another method, known as parametric MO, was proposed by \citeauthor{hugo2004}\cite{hugo2004}, which defines multiple discrete scenarios and solves the MO problem for each one, allowing for a comparison of the resulting Pareto sets. \citeauthor{tock2015}\cite{tock2015} showcased a probabilistic approach using Monte-Carlo simulations to quantify the probability distribution of each Pareto set, enhancing the uncertainty analysis within MO frameworks. For additional examples of different types of robustness in a multi-criteria context, we refer the interested reader to \citeauthor{ide2016}\cite{ide2016}.

Our new method presented here fit within the framework of adaptive scalarization-based algorithms for solving MO problems, as described in the context of chemical engineering in \citeauthor{bortz2014}\cite{bortz2014}. Furthermore, we interpret the robustification of the scalarization problems in a strict worst-case sense: the goal is to optimize the values of the objectives under the worst-case scenario, while feasibility of the solution must hold under any scenario. \replaced[id=R1]{In contrast to static worst-case robust optimization used in \citeauthor{bortz2014}\cite{bortz2014}, we now take the adaptability of the WSV into account.}{However, in contrast to standard worst-case robust optimization, we take the adaptability of the WSV into account.}

The paper is organized as follows: The next section provides a brief overview of how multi-criteria optimization problems under uncertainty can be appropriately reformulated and solved when there are HNV as well as WSV and the uncertainties realize in between. For solvability and performance reasons, the adaptive choice of scenarios from a reference discretization of the uncertainty set is proposed. Section 3 describes how to obtain the price of robustness, both for HNV and WSV.
%Robustness as additional objective is dealt with in section 4.
The conceptual findings are illustrated in Section 4 with an industrially relevant flow sheet. The paper ends with a conclusion and an outlook on future research.

\section{Multi-criteria, adjustable robust process design optimization}
\label{sec:MARO}

Tackling the task of choosing design and operating conditions as multi-criteria optimization (MO) problem has been shown to be very effective. For a detailed overview, we refer to Refs.~\citenum{clark1983,achenie2005,bortz2014}. By incorporating MO into an interactive decision support system, it is possible to explore the Pareto set and identify Pareto solutions that single-objective methods might miss \cite{bortz2014}. 

An MO problem is expressed as
\begin{equation}\label{eq:MOproblem}
	\min_{\replaced[id=R1]{(x,y) \in Z}{}} F(x,y,u_\text{nom})
\end{equation}

where $x$ (and $y$) denotes the HNV (WSV) along with their \replaced[id=R1]{common feasible domain $Z$}{respective feasible domains $X$ (and $Y$)}, \replaced[id=R1]{which is}{These feasible domains can be} specified explicitly or implicitly through inequality or equality constraints that are not detailed here. The function vector $F(x,y,u_\text{nom})=(f_1 (x,y,u_\text{nom} ),\dots,f_M (x,y,u_\text{nom}))^T$ includes $M$ objective functions. If the parameters $u$, which are included in the modeling and are often the source of uncertainty, e.g., estimated physical properties or raw material \replaced[id=R2]{prices}{prizes}, are set to some \emph{nominal} value $u_\text{nom}$, the outcome of the MO problem \eqref{eq:MOproblem} is the nominal Pareto set, which has significant importance for decision support in process design and optimization\cite{bortz2014,burger2014}.

\begin{figure}[h!]
    \centering
    % First subfigure
    \begin{subfigure}[b]{0.4\textwidth}
        \centering
        \includegraphics[width=\textwidth]{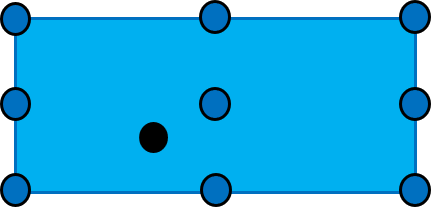}
        %\includegraphics[width=.35\textwidth]{gfx/BoxShapedRefDiscretization.png}
        %\hspace{1cm}
	%\includegraphics[width=.35\textwidth]{gfx/EllipticalRefDiscretization.png}
        \caption{}
        \label{fig:subfig1}
    \end{subfigure}
    \hfill
    % Second subfigure
    \begin{subfigure}[b]{0.4\textwidth}
        \centering
        \includegraphics[width=\textwidth]{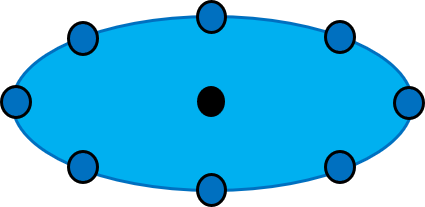}
        \caption{}
        \label{fig:subfig2}
    \end{subfigure}
    \caption{Reference discretizations for (a) a box-shaped and (b) an elliptical uncertainty set $U$: box vertices and face mids for the box-shaped uncertainty set and piercing points in coordinate directions as well as in direction of the space diagonals for the elliptical uncertainty set.
    The nominal scenario is marked in black and can, but does not have to lie in the center of the uncertainty set.}
    \label{fig:discretizations}
\end{figure}

To find Pareto optimal solutions hedging against uncertainties in the model parameters one realizes a \emph{multi-criteria adjustable robust optimization (MARO)} problem: 
\begin{equation}\label{eq:MARO}
	\min_{x \in \mathbf{R}^m} \max_{u \in U} \min_{y \in \mathbf{R}^n} F(x,y,u)\, \added[id=R1]{\text{ s.t. } g_l(x,y,u) \leq 0, l = 1,...,L} ,
\end{equation}
where the non-adjustability of the HNV $x$ and the adjustability of the WSV $y$ after realization of the uncertainties is taken into account. \added[]{Here, we restrict ourselves to fixed ($x$-independent) uncertainty sets $U$ and inequality constraints. Extensions to equality-constrained problems, which are common in process engineering, will be discussed later.} This formulation is a challenging tri-level optimization problem. In order to tackle it within common flowsheet simulators, it should be reformulated as an one-level optimization problem. To this end, the continuous uncertainty set $U$ is replaced by a discrete set of $K$ scenarios, the so-called \emph{reference discretization},
\begin{equation}\label{eq:discUset}
    \dot{U}_{\textrm{ref}} = \{ u_k \in U, k=1,\dots,K \},
\end{equation}
including the nominal scenario $u_{\text{nom}}$. To model the adjustability of the WSV $y$, a \emph{replication} is done: For each scenario $u_k$, an optimization variable $y_k$ is introduced. This replication and the consideration of finitely many scenarios enable to reformulate problem~\eqref{eq:MARO} as follows:
\added[id=R1]{ 
\begin{equation}\label{eq:discMARO}
    \min_{x\in \mathbf{R}^m, y_{1,\dots,K}\in \mathbf{R}^{K \cdot n}}\tilde{F}(x,y_{1,\dots,K}) \mathrm{ s.t. } \tilde{G}_l(x,y_{1,\dots,K}) \leq 0, l = 1,...,L
\end{equation}
}
with
\begin{equation}
    \tilde{F}(x,y_{1,\dots,K}) = \left(\tilde{f_1}(x,y_{1,\dots,K}),\dots,\tilde{f_M}(x,y_{1,\dots,K}) \right) ^T,
\end{equation}
\begin{equation}\label{eq:maxObjs}
    \tilde{f_j}(x,y_{1,\dots,K}) = \max_{k=1,\dots,K} \{f_j(x,y_k,u_k)\}, j = 1,...,M.
\end{equation}
\added[id=R1]{and
\begin{equation}
    \tilde{G}_l(x,y_{1,\dots,K}) = ( g_l(x,y_1,u_1),\dots,g_l(x,y_K,u_K) )^T, l=1,...,L.
\end{equation}
}
Problem~\eqref{eq:discMARO} is merely an approximation because presumably the worst-case points of the continuous uncertainty set $U$ obtained in Eq.~\eqref{eq:MARO} are not the same as those resulting from the reference discretization $\dot{U}_{\textrm{ref}}$ considered in Eq.~\eqref{eq:discMARO}. However, the discretization $\dot{U}_\textrm{ref}$ of $U$ makes it possible to compute an approximate solution of Problem~\eqref{eq:MARO} by enumerated optimization instead of applying heavy-weighted methods of global optimization to solve the $2^{nd}$- and $3^{rd}$-level problems in Problem~\eqref{eq:MARO}. Examples for such reference discretizations are shown in Fig.~\ref{fig:discretizations}. In order to keep the approximation error small we suggest a reference discretization $\dot{U}_\textrm{ref}$ which is as fine as computational possible. 
\added[id=R2]{The finer the reference discretization, the more accurate the results, but also the more computationally intensive the calculation. We therefore suggest first calculating with a coarse discretization, e.g., the corner points and the center point of an uncertainty box, to get a first approximation of the worst-case solution, and then calculating with a refined discretization — possibly only locally around the worst-case scenarios — to check and improve the accuracy of the approximation.}
\added[id=R12]{As a result of using finitely many scenarios, it is also possible to consider non-polyhedral, e.g., elliptical, or completely arbitrary uncertainty sets — applying the same solution approach introduced next.}

The \emph{replication} or \emph{function value approach} is grounded on the assumption that the WSV $y$ can be adjusted to previously realized uncertainties $u$, i.e., on the (functional) dependence of $y$ on $u$. However, since the functional relationship is generally not known, the replication approach takes the function values $y_k = y(u_k)$ into account instead of considering the function $y(\cdot)$ itself. Another option would be to use a linear, piecewise linear, or quadratic ansatz for the unknown function $y(\cdot)$. This is referred in the literature as \emph{decision rule approach} and can be explored in more depth, e.g., in Refs.~\citenum{yanikouglu2019,schneider2024} and the references therein.

\subsection{Adaptive scenario choice in MARO}
\label{sec:adaptive_scenarios}

The straightforward way to solve problem~\refeq{eq:discMARO} is to consider all scenarios $u_k, k=1,\ldots,K,$ of the reference discretization at once, leading to large-scale optimization problems. In general, however, not all reference scenarios are relevant for determining a worst-case solution; often only a small subset of them is needed. Therefore, we choose the scenarios adaptively. This procedure is inspired by the Blankenship-Falk algorithm \cite{blankenship1976} and works as follows:

\subsubsection{Algorithm~1: Adaptive scenario choice for worst-case Pareto point computation}

\paragraph{Step~(0)}\textit{Initialization:} Choose an initial subset $\dot{U}_0 \subseteq \dot{U}_{\text{ref}}$, e.g., $\dot{U}_0 = \l\{ u_{\text{nom}} \r\}$, and apply a scalarization technique with scalarization parameters $\alpha$ to transform the discretized MARO problem~\refeq{eq:discMARO} to a single-objective problem with objective function $f_\alpha$ (see, e.g., Ref.~\citenum{bortz2014}, for possible scalarizations). Set the iteration index $i:=0$.

\paragraph{Step~(1)}\textit{Computation of (approximate) worst-case solution:} Solve the resulting scalarized and discretized MARO problem yielding an optimal solution $\l( x^*_{(i)}, y^*_{(i)} \r)$ and optimal function values $f^*_{j,(i)}$.

%\color{red}

\paragraph{Step~(2)}\textit{Scenario re-optimization:} For $x^*_{(i)}$ fixed and each scenario $u_k, k = 1,...,K,$ of the reference discretization solve the (scalarized) scenario re-optimization problem:
\begin{equation}\label{eq:ScenReOpt}
    \min_{y \in \mathbf{R}^n} f_{\alpha}(x^*_{(i)}, y, u_k) \textrm{ s.t. } g_l(x^*_{(i)},y,u_k) \leq 0, l = 1,...,L,
\end{equation}
yielding on optimal solution $\l( x^*_{(i)}, y^*_{(i),k} \r)$.

\paragraph{Step~(3)}\textit{Determination of worst-case scenarios:}

\textbf{(a)} For each objective $f_j, j = 1,...,M$: If the worst objective value on the reference discretization is worse than the current worst-case objective value, i.e.,
\begin{equation}\label{eq:curObjWCscen}
    \max_{k=1,...,K} f_j(x^*_{(i)}, y^*_{(i),k}, u_k) > f^*_{j,(i)},
\end{equation}
add objective worst-case scenario $u^*_{(i,f_j)} = \mathrm{argmax}_{k=1,...,K} f_j(x^*_{(i)}, y^*_{(i),k}, u_k)$ to the current set of worst-case scenarios $\dot{U}_{i+1} = \dot{U}_i \cup \{ u^*_{(i,f_j)} \}$.

\textbf{(b)} For each constraint $g_l, l = 1,...,L$: If the worst constraint value on the reference discretization is positive, i.e.,
\begin{equation}\label{eq:curConsWCscen}
    \max_{k=1,...,K} g_l(x^*_{(i)}, y^*_{(i),k}, u_k) > 0,
\end{equation}
add constraint worst-case scenario $u^*_{(i,g_l)} = \mathrm{argmax}_{k=1,...,K} g_l(x^*_{(i)}, y^*_{(i),k}, u_k)$ to the current set of worst-case scenarios $\dot{U}_{i+1} = \dot{U}_i \cup \{ u^*_{(i,g_l)} \}$.

\noindent Replace $i$ by $i+1$. If there are new worst-case scenarios, i.e., $\dot{U}_{i+1} \supsetneq \dot{U}_{i}$, go back to Step~(1), otherwise an optimal solution $\l( x^*, y^* \r) = \l( x^*_{(i)}, y^*_{(i)} \r)$ of \eqref{eq:discMARO} w.r.t. to the reference discretization $\dot{U}_{\textrm{ref}}$ is found.

%\color{black}

%\item[vi)] Go back to i) and change the scalarization parameters $\alpha$. These are chosen adaptively to obtain the Pareto boundary with a desired accuracy \citenum{bortz2014}. Iterate until this accuracy goal has been achieved.

\deleted[id=R1]{If the process design optimization task contains inequality constraints, which depend on the uncertainties, the following additions have to be made:}
\deleted[id=R1]{
\begin{itemize}
    \item The inequality constraints must be replicated for each scenario of the set of current worst-case scenarios $\dot{U}_{i}$ in Step~(1),
    \item the worst-case scenarios regarding the constraints, i.e., those that violate the constraints most severely on the reference discretization, have to be determined in Step~(2), too, and added to the current set of worst-case scenarios.
\end{itemize}
}

Equality constraints can either be included in the flow sheet equation system or eliminated by additional process variables (see, e.g., \citeauthor{barton2015}\cite{barton2015} for details).

\added[id=R12]{The running time of the algorithm depends heavily on the size of the reference discretization:
\begin{itemize}
    \item In Step~(2), a generally nonlinear optimization problem must be solved for each scenario of the reference discretization.
    \item For each scenario included in the current discretization $\dot{U}_{(i)}$, a set of WSV and constraints for must be added to the discretized MARO problem \eqref{eq:discMARO}, which increases the size of the discretized MARO problem rapidly.
\end{itemize}
}

%\added[id=R1]{If one wants to consider continuous uncertainty sets, one would have to solve a two-stage optimization problem to global optimality in Step (2). Under certain assumptions, e.g., linearity in all variables $x$, $y$, and $u$, this would be more efficient as with CCG-algorithm\cite{Zeng2013}, but under general conditions it is not.}

\added[id=R1]{In the case where the HNV and WSV enter linearly into the problem formulation (the WSV being also mixed-integer) and the uncertainty set is polyhedral or finite, Alg.~1 corresponds to the column-and-constraint generation (C\&CG) algorithm for two-stage robust optimization\cite{Zeng2013}. Under these conditions, it would also be possible to consider the uncertainties continuously and calculate the exact worst-case solution. To do this, instead of Steps~(2) and (3) a two-stage linear optimization problem has to be solved or — assuming a fixed ($x$-independent) uncertainty polytope — its vertices would have to be considered as reference discretization in Steps~(2) and (3).}

\begin{figure}
	\includegraphics[width=.24\textwidth]{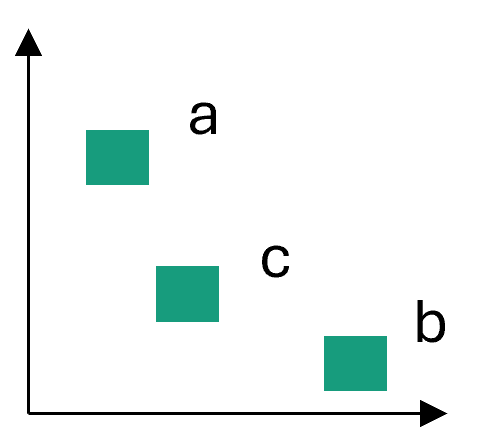}
	\includegraphics[width=.24\textwidth]{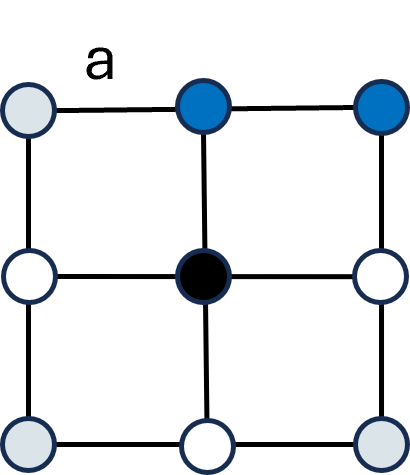}
	\includegraphics[width=.24\textwidth]{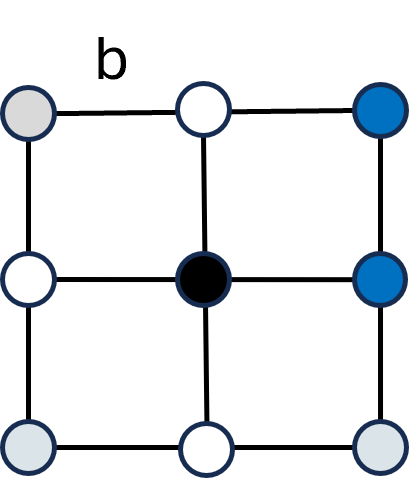}
	\includegraphics[width=.24\textwidth]{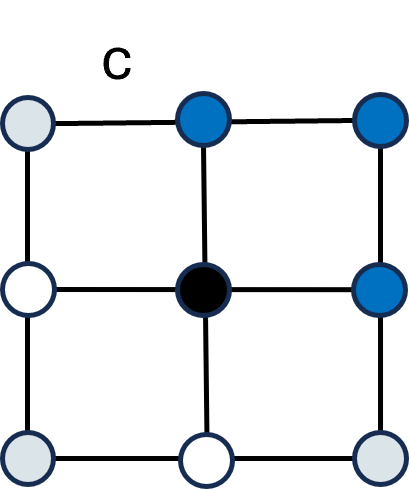}
	\caption{Sketch of adaptive scenario choice for three Pareto points a, b, c. The scenarios are chosen from a reference discretization (black, gray, and white points, where the black point represents the nominal scenario) according to Alg.~1. For Pareto point c the union of the (blue) worst-case scenarios for point a and point b is chosen as initial set of worst-case scenarios.}
    \label{fig:ParetoDisc}
\end{figure}

In addition to the adaptive determination of the worst-case Pareto solution, the worst-case Pareto front, i.e., the set of worst-case Pareto solutions, can also be determined adaptively \cite{bortz2014}. The adaptive determination of each worst-case Pareto solution does not have to \replaced[id=R2]{start with}{be based on} a nominal Pareto solution, i.e., with $\dot{U}_0 = \l\{ u_{\text{nom}} \r\}$. Instead, the \replaced[id=R2]{worst-case}{WC} scenarios determined so far can be used to initialize the set $\dot{U}_0$ (see Fig.~\ref{fig:ParetoDisc} for an illustration). This significantly speeds up the refinement process (alternating between Steps~(1) and (2)) and often leads to the refinement being terminated after just one alternation.

The adaptive determination of the worst-case Pareto solution can be carried out in the same way if there are no WSV $y$. In this \replaced[id=R1]{ static}{, standard} robust case, the $y$ are omitted as optimization variables and the optimization problem in Eq.~\ref{eq:ScenReOpt} is reduced to \replaced[id=R12]{an evaluation of the objective and constraint functions}{a function evaluation}. 

%For an adaptive determination of the worst-case solution in the single-criterion case with \highlight[id=NA]{regret objective function and affine decision rule} instead of replication we refer the reader to the work~\citenum{schneider2024}.
%\todo[inline]{NA: Ist das hier zu verstehen? Weder regret noch affine decision rule vorher erklärt.}

\section{Price of robustness}

The following section deals with quantifying the compromises with respect to the performance indicators, or objective functions, which have to be made in order to hedge against the uncertainty. \replaced[id=R1]{The concept we apply was introduced in \citeauthor{bertsimas2004}\cite{bertsimas2004} as the price of robustness.}{We call these compromises the price of robustness.}

If engineers formulate the process design problem as a worst-case robust optimization problem, for example, in the form of a MARO problem \eqref{eq:MARO} they decide to simultaneously hedge against all scenarios from a chosen scenario set $U$. However, by accounting for all scenarios that may happen, the performance of this solution in a particular, realized scenario will generally be worse than a solution tailored to this particular scenario only would have been.
%\todo[inline]{MW: Muss das ueberaupt erklaert werden, oder ist es eh klar, weil es um robuste Optimierung geht?}  

%%In a practical setting, uncertainties are often construed as deviations from a specific nominal value. The scenario in which all uncertain parameters adopt their nominal value is called the nominal scenario ($u_{nom}$). The nominal scenario may be seen as “more likely” to occur than a more extreme scenario at the edge of the uncertainty set. If this is the case,

Usually, the nominal scenario $u_{\text{nom}}$ is of particular relevance as uncertainties are varied around its values. Thus, it seems useful to consider the loss of performance in the nominal scenario incurred by hedging against all other scenarios. This loss of performance of a robust solution in the nominal scenario is called  \emph{price of robustness}.  

\subsection{Quantifying the price of robustness for MARO} 

In accordance with the
%multi-objective problem
formulations MO (\ref{eq:MOproblem}) and MARO (\ref{eq:MARO}), where the performance of a solution is conceived as an $m$-dimensional vector of performance indicators $F$, we also conceive the price of robustness as an $m$-dimensional vector $p^R$, where each entry denotes the loss in the corresponding performance indicator. 

We apply the adaptive scenario choice approach described in the previous section and define 
\begin{equation}
    F^{MARO}(x^*):= \tilde{F}(x^*,y^*_{1,\dots,K})
\end{equation}
as the image of an optimal solution $(x^*,y^*_{1,..,K})$ to the approximated reformulation of the MARO problem, as defined in (\ref{eq:discMARO}) and (\ref{eq:maxObjs}).    

For the MARO problem (\ref{eq:MARO}), hedging against all scenarios occurs in the choice of the HNV $x$, while the WSV $y$ can be chosen according to the realized scenario. With this consideration, the performance of the robust HNV solution  $x^*$ in the nominal scenario is the solution to the optimization problem
\added[id=R2]{ 
\begin{equation}\label{eq:NSRpt}
    \min_{y \in \mathbf{R}^n} F(x^*,y,u_\text{nom}) \mathrm{ s.t. } g_l(x^*,y,u_\text{nom}) \leq 0, l = 1,...,L
\end{equation} 
}
where $y$ is optimized both in response to the choice of $x^*$ and to the realization of the nominal scenario. We call \refeq{eq:NSRpt} the \emph{nominal scenario re-optimization (NSR) problem}.

The NSR problem is an MO problem itself, and has - in general - an infinite number of solutions. As a direct consequence of the MARO problem definition \eqref{eq:MARO}, any solution $F^{MARO}(x^*)$ for $x^{*}$ is an upper bound for what is achievable in the NSR. Generally, since in Eq.~\eqref{eq:NSRpt} $u_{\text{nom}}$ is taken instead of the worst-case as in \eqref{eq:MARO}, a better performing WSV $y$ can be found. When deciding on the WSV $y$ in the NSR, the engineer will make a particular choice $y^{NSR}$, depending on the relative value he assigns to the different performance indicators (\replaced[id=R2]{their}{her} preference). We denote the resulting particular preferred Pareto point as
\begin{equation}
\label{eq:NSRpt_point}
F^{NSR}(x^*):=F\l(x^*, y^{NSR}, u_{\text{nom}}\r). 
\end{equation}
Compared to the solution $F^{MARO}(x^*)$, the solution to the re-optimization problem $F^{NSR}(x^*)$ with fixed HNV $x^*$ is improved by
\begin{equation}
d:=F^{NSR}(x^*)-F^{MARO}(x^*). 
\end{equation}

\begin{figure}
	\includegraphics[width=.5\textwidth]{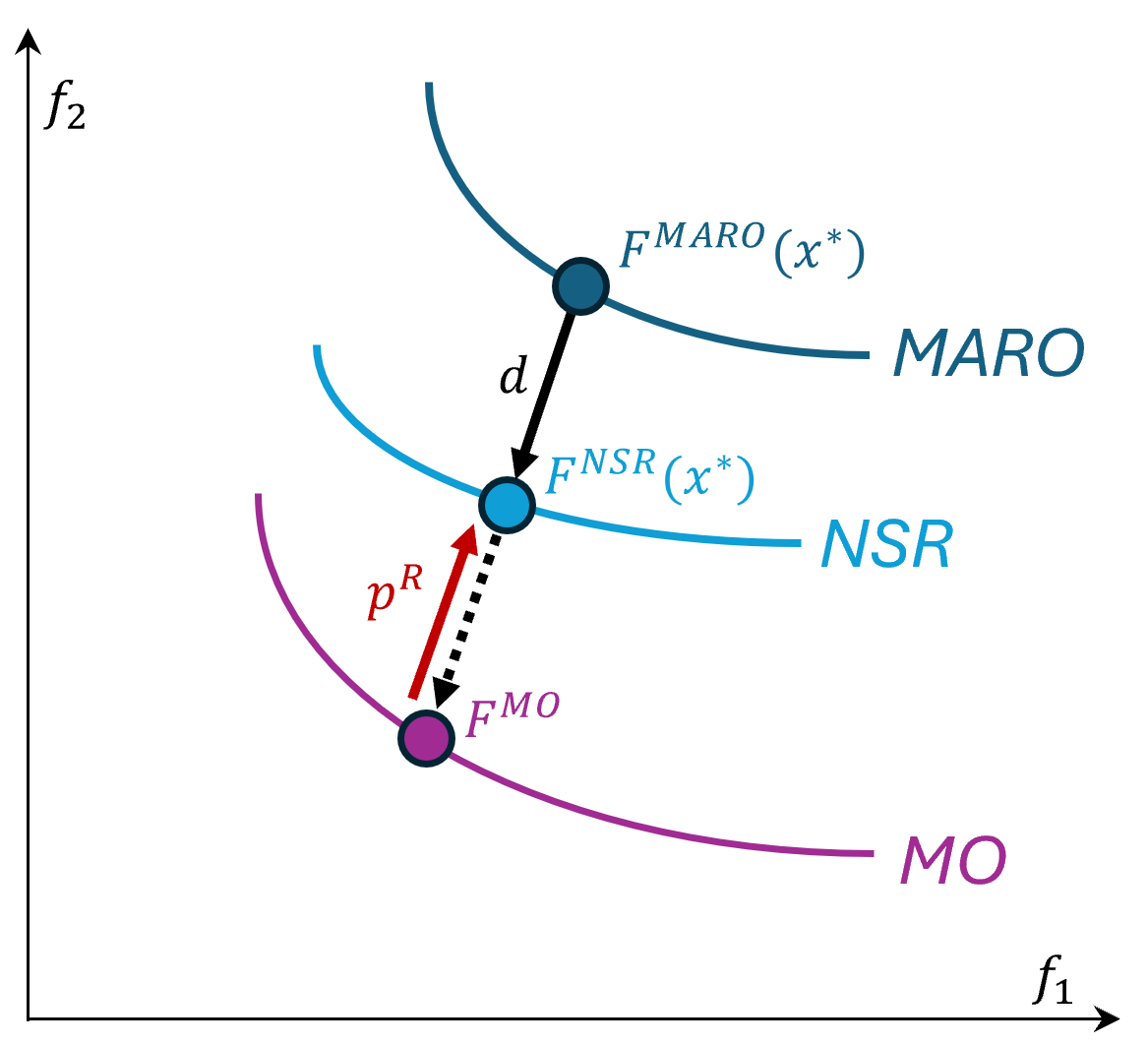}
	\caption{Illustration of the price of robustness in the MARO setting. The image of a MARO optimized HNV $x^*$ is a point $F^{MARO}(x^*)$ on the MARO Pareto front. Re-optimization of WSV $y$ for the nominal scenario $u_{\text{nom}}$ leads to an improved point $F^{NSR}(x^*)$. The direction of improvement can be extended to a point $F^{MO}$ on the nominal front. The price of robustness $p^R$ is the distance from $F^{MO}$ to $F^{NSR}(x^*)$.}
\label{fig:por}
\end{figure}

As $F^{NSR}(x^*)$ was chosen according to the preference of the engineer, we can assume that $d$ represents the preferred relative distribution of improvements among the performance indicators, compared to $F^{MARO}(x^*)$. Analogously, we can find a Pareto efficient solution $(x^{MO}, y^{MO})$ to the non-robust problem, represented by a point $F^{MO}=F(x^{MO}, y^{MO})$ on the nominal Pareto front, that improves the performance indicators at the same ratio, again compared to $F^{MARO}(x^*)$. The geometric interpretation is that the intersection of the ray starting from $F^{MARO}(x^*)$ in the direction $d$ meets the nominal Pareto front of problem \eqref{eq:MOproblem} in point $F^{MO}$:
%\todo[inline]{MW: Man muss ich hier viel denken. Evtl. koennte man die Definitionen von $F_{MARO}$, $F_{NSR}$ nochmal ausfuehrlicher als Loesungen des jeweiligen Optimierungsaufgabe hinschreiben, und deutlich machen, dass man das Design $x^*$ aus MARO nimmt, mihilfe des NSR eine Richtung berechnet, und dass dann aber der Schnittpunkt mit der gewoehnlichen/Nominal-Paretofront $F^*$ i.A. ein anderes Design (und Betriebsvariablen) hat.}

\begin{equation}
\label{eq:comp_point}
\begin{split}
    F^{MO} & := F^{MARO}(x^*) + \alpha^* d.
\end{split}
\end{equation} 
As for the point $F^{MO}$ on the nominal front both HNV $x^{MO}$ and WSV $y^{MO}$ optimally account for the nominal scenario, we can expect $\alpha^* > 1$ in equation (\ref{eq:comp_point}), meaning that a larger improvement of $F^{MO}$ over $F^{MARO}(x^*)$ can be attained than the improvement $d$ from the re-optimization of the WSV alone. With the above definitions, the vector-valued price of robustness is the loss incurred by robustifying the HNV:
\begin{equation}
\label{eq:por}
p^R = F^{NSR}(x^*) - F^{MO}.
\end{equation}
Fig.~\ref{fig:por} shows a schematic illustration of the price of robustness for the MARO problem.

%\todo[inline]{MW: In der Grafik ist der Vektor d falsch herum eingezeichnet. Er zeigt von NSR nach MARO!}
%\todo[inline]{NA: Achsen bis zum Bildende faende ich schoener.}

\subsection{Real-time approximation and interactive display of the price of robustness} \label{sec:interactive}

Using a decision support system, 
%integrated into Chemasim,
engineers can explore the different Pareto optimal solutions to the MARO problem. This navigation of the Pareto front is steered through a panel of sliders\cite{bortz2014}, each of which facilitates the targeted improvement and restriction of one of the performance indicators. 

%Interactive navigation is facilitated by pre-calculation and subsequent linear interpolation of a well distributed set of solutions. Pareto front approximation algorithms ensure that the pre-calculated solutions are well distributed and the interpolation error is small \ref{}.

The price of robustness can be incorporated into this framework by, say for performance indicator $f_k$, visualizing both $f^{NSR}_k(x^*)$ and $f^{MO}_k$ as markers on the corresponding slider. Then, the distance between these markers indicates the price of robustness for the $k$-th performance indicator. While exploring the Pareto front, both markers update continuously on the sliders when $x^*$ is changed and the price of robustness can be observed in real time.

Smooth real-time manipulation of the navigated point $F^{MARO}$ is facilitated by linearly interpolating between pre-calculated solutions that approximate the Pareto front of the MARO problem \eqref{eq:MARO} (see \citeauthor{monz2008}\cite{monz2008} for details). The real-time update of the markers works very similarly. When computing the approximating solutions, the re-optimization problems are solved upfront. During navigation, the interpolation coefficients defining the navigated point are also used to obtain $F^{NSR}(x^*)$. 

\replaced[id=R2]{Similar to the navigated point $F^{MARO}$, the non-robust point $F^{MO}$ (definition \eqref{eq:comp_point}) is approximated by finding a point on the nominal front approximation rather than the nominal front itself. The exact algorithm to find the approximated point $\hat{F}^{MO}$ depends on whether the nominal front is convex or not. If the nominal front is convex, its approximation is the convex hull of the pre-calculated points, and the convex coefficients that define $\hat{F}^{MO}$  can be obtained by solving a linear program (very similar to the simultaneous navigation on multiple Pareto fronts described in Ref.~\citenum{patchnav2021}). If the nominal front is non-convex, its approximation is instead obtained from a Delaunay triangulation of the pre-calculated points (see \citeauthor{nowak2022}\cite{nowak2022} for details), and a bisection algorithm is used to search for $\hat{F}^{MO}$ as the point where the ray emanating from $F^{MARO}$ in direction $d$ intersects with a face of the Pareto front approximation. Irrespective of the different algorithmic approaches for the convex and non-convex case, $\hat{F}^{MO}$ can be calculated fast enough to facilitate a real time update of the price of robustness during navigation.}{For the calculation of the non-robust point $F^{MO}$, definition \eqref{eq:comp_point} is translated to an optimization problem to find $F^{MO}$ as a convex combination of neighbored pre-calculated nominal Pareto solutions, which can be solved in real time (very similar to the simultaneous navigation on multiple Pareto fronts described in Ref.~\citenum{patchnav2021}).}

\begin{figure}
	\includegraphics[width=.75\textwidth]{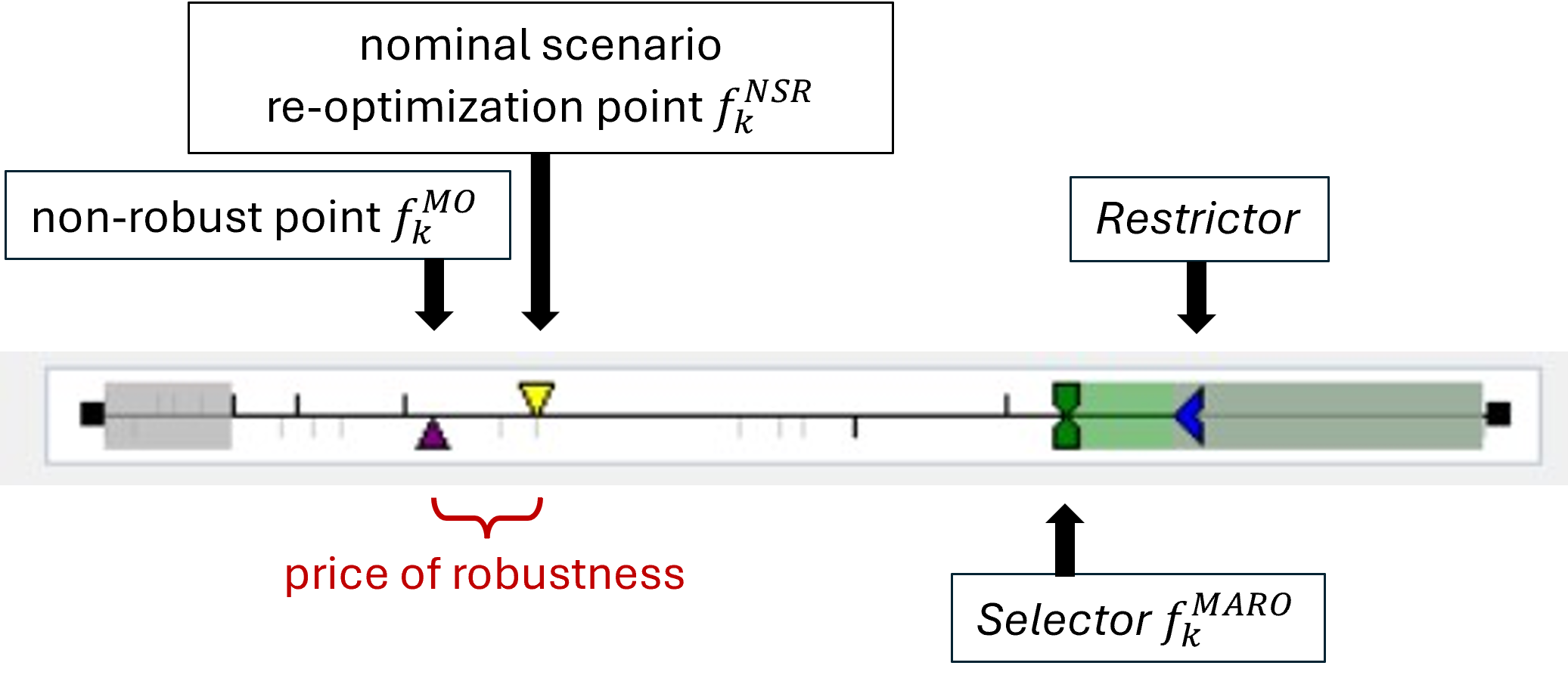}
	\caption{Illustration of the price of robustness displayed on one of the navigation sliders. The worst-case objective value for the MARO solution is represented by the green hourglass selector, and can be improved by dragging the selector to the left. Yellow and purple markers indicate the objective value for the re-optimized solution (NSR) and the non-robust solution (MO), respectively. Their difference is the price of robustness for this particular objective.}
\label{fig:sliders_schem}
\end{figure}

\section{\replaced[id=R12]{Case study}{A numerical example}: Separation of Methanol and Methyl formate}

We now demonstrate the adaptive calculation of worst-case Pareto solutions and the computation as well as the visualization of the price of robustness using a small flow sheet example.

\subsection{Setup}

\begin{figure}[h]
	\includegraphics[width=.55\textwidth]{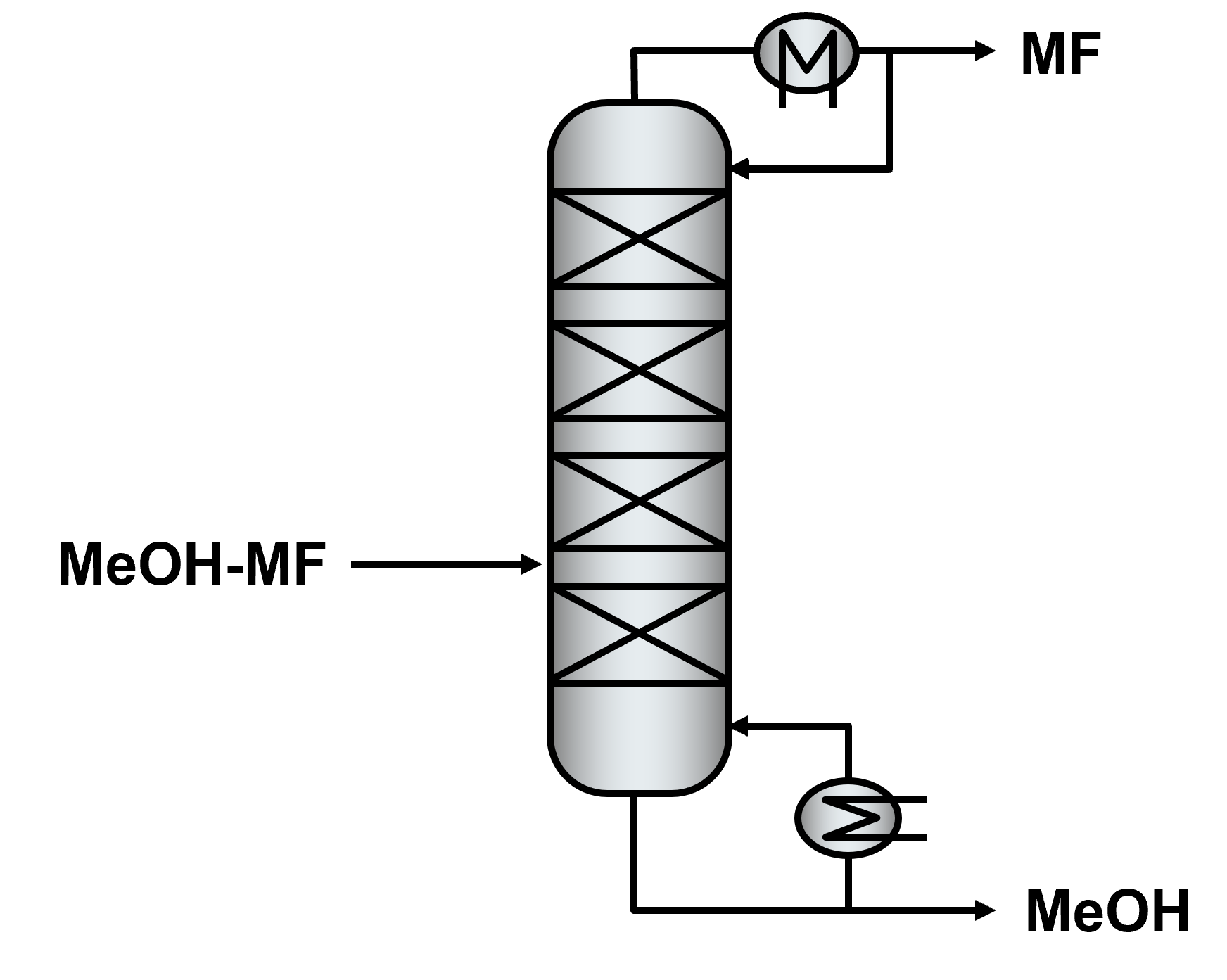}
	\caption{Separation of a Methanol - Methyl formate mixture using a distillation column.}
    \label{fig:MeOHMFdistColumn}
\end{figure}

A standard distillation column
%with a reboiler at the bottom and a condenser at the top and with bottom and top product streams
is considered to separate a binary zeotropic mixture of \emph{Methanol (MeOH)} and \emph{Methyl formate (MF)}. The light boiler MF is withdrawn at the top and the heavy boiler MeOH at the bottom of the column (see Fig.~\ref{fig:MeOHMFdistColumn}).

The challenge consists in identifying HNV and WSV in order to find the Pareto set for minimizing the \emph{capital expenditures (CAPEX)} and the \emph{operating expenditures (OPEX)} (per ton), assuming uncertainties in the load of the column, the feed composition and the activity coefficients. \added[id=R1]{In the case of these two cost objectives, we've decided not to aggregate them a priori into one overall cost function, because the multi-criteria setting is a key feature of our solution approach and the decision-maker should be able to quantify the trade-off between the two costs in order to assess changes in the costs.}

As HNV $x$ (i.e., the variables related to the design of the separation process) we choose the number of equilibrium stages $N$, the feed stage $N_f$, the diameter $D$ of the column, and the areas of the heat exchangers in the reboiler and in the condenser, $A_r$ resp. $A_c$. The considered WSV $y$ (i.e., the operational variables of the column) are the reflux ratio $R_V$, and the specific heat duty of the reboiler $\dot{Q}_r$. The (uncertain) load $l$ of the column is defined as feed mass flow rate divided by 8000 kg/h; the uncertainty of the feed composition is captured by the mass fraction $w_{MF}$ of MF in the feed and the uncertainty in the activity coefficients is modeled by the thermodynamic factor approach described in Ref.~\citenum{burger2017}. For the units and ranges of the HNV, the WSV, and the uncertain quantities, see Tab.~\ref{tab:uncertainties}. 

\begin{table}
    \centering
    \begin{tabular}{|l|c|c|c|}
        \hline
        \textbf{Here-\&-now variables} & \textbf{unit} & \textbf{range} & \textbf{initial value} \\
        \hline
        Number of equilibrium stages $N$ & - & 10 - 150 & 33 \\
        Feed stage $N_f$ & - & 3 - 40 & 5 \\
        Column diameter $D$ & m & 0.8 - 2.0 & 1.09 \\
        Heat exchanger area in the reboiler $A_r$ & m\textsuperscript{2} & 50 - 1000 & 216.98 \\
        Heat exchanger area in the condenser $A_c$ & m\textsuperscript{2} & 50 - 1000 & 191.91 \\
        \hline
        \textbf{Wait-\&-see variables} & & & \\
        \hline
        Reflux ratio $R_V$ & - & 0.5 - 2.0 & 0.74 \\
        Specific heat duty of the reboiler $\dot{Q}_r$ & kWh/kg & 0.0625 - 0.375 & 0.21 \\
        \hline
        \textbf{Uncertain (model) parameters} & & & \textbf{nominal value} \\
        \hline
         Load $l$ & - & 0.6 - 1.2 & 1.0 \\
         Mass fraction of MF in feed $w_{MF}$ & kg/kg & 0.78 - 0.82 & 0.8 \\
         Thermodyn. factor $F_{12}$ for act. coefficients & - & 0.9 - 1.1 & 1.0 \\
         \hline
    \end{tabular}
    \caption{HNV, WSV, and uncertain quantities with units, ranges and initial/nominal values for the MeOH-MF separation problem.}
    \label{tab:uncertainties}
\end{table}

%\begin{figure}
%	\includegraphics[width=.45\textwidth]{gfx/fig_scenarios.jpg}
%	\caption{Placeholder: Scenario set for MeOHMF.}
%\end{figure}

\noindent Furthermore, the following restrictions should hold for the selected HNV and adjusted WSV in each scenario:
\begin{itemize}
    \item a (minimum) purity requirement each for the products MeOH and MF, $w_{MeOH,bot}^{min}$ resp. $w_{MF,top}^{min}$
    \item an upper bound for the heating resp. cooling capacity of the reboiler and the condenser, $Q_r^{max}$ resp. $Q_c^{max}$, and
    \item a restriction on the F-factors of the distillation column $F^{max}$ (not to be confused with $F_{12}$ used in Tab.~\ref{tab:uncertainties} to model the uncertain activity coefficient, Ref.~\citenum{burger2017})
\end{itemize}

\added[id=R12]{The distillation column with reboiler and condenser was modeled in CHEMASIM, BASF's in-house flowsheet simulator, using equilibrium stages, i.e., MESH equations. These are inherently nonlinear. On top of that we have BASF's confidential cost functions. The constraints also depend generally nonlinearly on the variables under consideration. Thus, for all HNV, WSV, and uncertainties, there are nonlinear relationships in the objective and constraint functions. This means that even for box uncertainties, a discretization of the uncertainty box is necessary. Furthermore, we will consider the case of an elliptical uncertainty set.}

\subsection{Calculation of the MARO front using the adaptive scheme}

Using the adaptive scenario selection algorithm Alg.~1 we calculate the Pareto front approximation of the resulting MARO problem \added[id=R12]{for the MeOH-MF separation introduced above. For an in depth investigation of the results, we employ one particular discretization scheme, namely a box‐shaped uncertainty set comprising its vertices and face midpoints. At the end of this section, we consider additional discretization schemes to show the substantial reduction in computation time achieved by the adaptive scenario selection algorithm.}     

%\todo[inline, backgroundcolor=blue!20!white]{Pareto front approximation}

From our previous work \cite{schwientek2025} we know that the Pareto fronts of optimal compromises between CAPEX and OPEX for the MeOH-MF separation problem are convex, both for the nominal, non-robust problem formulation \eqref{eq:MOproblem} as well as for the MARO formulation \eqref{eq:discMARO}. This allows us to employ the Sandwiching algorithm \cite{Serna2012,Bokrantz2012}, to obtain a set of representative Pareto optimal solutions that approximate the Pareto fronts. In the Sandwiching algorithm, each calculated Pareto solution is the solution to a weighted-sum scalarization problem. \added[id=R12]{However, the choice of scalarization is not critical; the Pareto front of the MARO problem could also be calculated using other scalarization methods such as Epsilon-Constraint or Pascoletti-Serafini, and in the case of a non-convex problem, those would be even more appropriate.} 
%\todo{MW: Ist das a-priori klar?}
%\todo{NA: Ref. ergaenzen}
%\todo[inline]{PasS: Sounds odd this way. Would either say that we know from previous work that it is convex or that we use the hybrid stuff and as it is convex, this is simply the SW.}

%\todo[inline, backgroundcolor=blue!20!white]{Initial reference discretization}

As reference discretization, we employ the box-shaped version schematically shown in Fig.~\ref{fig:discretizations}(a), but in a 3-dimensional setting (see Fig.~\ref{fig:MeOHMFscenarios}). As we combine the minimal, maximal, and mid value for each of the three uncertain parameters listed in Tab.~\ref{tab:uncertainties},
%\deleted[id=R12]{
and additionally include the nominal scenario $u_{nom} = (1.0, 0.8, 1.0)$,
%}
we obtain
%\replaced[id=R12]{$3^3 = 27$}{
$3^3+1 = 28$
%}
scenarios in total.
%\added[id=R12]{, where the box center $(0.9, 0.8, 1.0)$ is taken as nominal scenario.}
%\todo[inline]{PasS: Would add a note here again to make things clear, i.e., that we take min, mid and max for each param and then any combination giving the $3^3$. The nominal scenario is always considered even though it is here the combination of all mid points and thus already included.}

%\todo[inline, backgroundcolor=blue!20!white]{Application of adaptive scheme and obtained worst case scenarios}

For the calculation of the \replaced[id=R2]{worst-case}{WC} extreme compromises, Alg.~1 was initialized with only the nominal scenario as starting scenario set, i.e., nominal extreme compromises were computed initially. After 1 refinement of the initial scenario set $\dot{U}_0 = \{ u_{nom} \}$ (based on the reference discretization with 28 scenarios), the worst-case extreme compromises were determined using 5 scenarios 
%(3, 21, 25, 27, 28/N)
for CAPEX and 4 scenarios 
%(6, 25, 27, 28/N)
for OPEX.

%\begin{figure}[H]
%    \centering
%    \includegraphics[width=.5\textwidth]{gfx/UncBox_MeOH-MF.png}
%        \captionof{figure}{Reference discretization used for the MeOH-MF example. All scenarios that appear as worst-case in one of the scalarizations are marked in red.}
%        \label{fig:MeOHMFscenarios}
%\end{figure}

The subsequent weighted-sum (WS) scalarizations were initialized with the union of the \replaced[id=R2]{worst-case}{WC} scenario sets of the extreme compromises, six in total (incl. the nominal scenario). \replaced[id=R12]{Three}{6} WS runs were performed to achieve the desired approximation quality. In these runs, the set of \replaced[id=R2]{worst-case}{WC} scenarios did not have to be refined in order to calculate \replaced[id=R2]{worst-case}{WC} Pareto solutions.

%Overall, the worst-case Pareto front could be approximated with few scenarios (4-6 per Pareto point; $\sim20\%$ of reference scenarios) and one refinement at maximum.
Overall, the worst-case Pareto front could be approximated with only a small subset ($\sim 20\%$) of the reference scenarios, which allows for major reductions in the computational time at vanishing losses in the approximation quality (see below for more details).

Fig. \ref{fig:MeOHMFscenarios} and Tab.~\ref{tab:scenlist} show which scenarios were finally identified as worst-case scenarios among the scenarios from the reference discretization.

%At each Pareto point, scenario 25 represents the WC scenario for the operating costs, while scenario 27 represents the WC scenario for the restrictions $Q^{max}_c$, $Q^{max}_r$, and $F^{max}$. For the purity requirements $w_{MeOH,bot}^{min}$ and $w_{MF,top}^{min}$, the WC scenario varies between scenarios 3, 6, 21, and 27 across the Pareto points. The investment costs do not have a worst-case, because they do not depend on the uncertainties.
The investment costs do obviously not have a worst-case, because they only depend on the HNV and thus do not change if adjusting the WSV according to the realization of the uncertainty. At each Pareto point, scenario 25 represents the \replaced[id=R2]{worst-case}{WC} scenario for the operating costs, as high thermodynamic factors and high feed fractions of MF represent the limiting case for the operating costs due to a substantially shrinking (and almost vanishing) width of the two-phase-region at high MF mole fractions (note however that the load does not have an influence on OPEX due to its definition per ton of product so that there are equally bad scenarios). Following this discussion, scenario 27 being the unique \replaced[id=R2]{worst-case}{WC} for the operating range of the equipment (represented by $Q^{max}_c$, $Q^{max}_r$, and $F^{max}$) is intuitive, as it further considers the high load case while having the most unfavourable thermodynamic factor and MF feed fraction in terms of operating costs (that directly stem from high energy input and thus high refluxes and vapour loads). For the purity requirements $w_{MeOH,bot}^{min}$ and $w_{MF,top}^{min}$, the \replaced[id=R2]{worst-case}{WC} scenario varies between scenarios \replaced[id=R12]{12 and 21}{3, 6, 21, 27} across the Pareto points. This loosely matches the engineering intuition (e.g., from simple consideration of a McCabe-Thiele diagram) that the MF content has a monotonic impact on achievable purities at top and bottom. However, as we simultaneously consider uncertainties in the VLE-behavior (with the tendency that higher thermodynamic factors shrink the two-phase-region at high MF mole fractions but widen it at low fractions), there is no clear \replaced[id=R2]{worst-case}{WC} for the purity requirements.

\definecolor{fghblue}{RGB}{0,91,127}
\definecolor{myred}{RGB}{247,15,34}
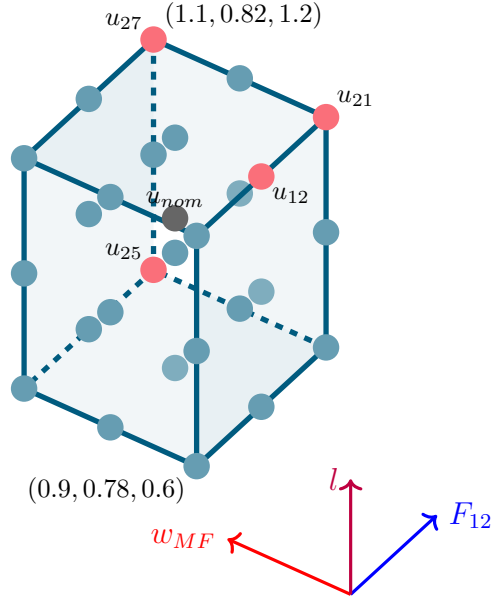
\begin{figure}[H]
    \centering
  \begin{tikzpicture}[tdplot_main_coords]

    % Dimensions
    \def\W{3}
    \def\D{2.8}
    \def\H{4}

    % Corner coordinates
    \coordinate (A) at (0,0,0);
    \coordinate (B) at (\W,0,0);
    \coordinate (C) at (\W,\D,0);
    \coordinate (Dpt) at (0,\D,0);

    \coordinate (E) at (0,0,\H);
    \coordinate (F) at (\W,0,\H);
    \coordinate (G) at (\W,\D,\H);
    \coordinate (Hpt) at (0,\D,\H);

    % Face centers as little balls
    \coordinate (BottomC) at ($ (A)!.5!(C) $); % bottom face center
    \coordinate (TopC)    at ($ (E)!.5!(G) $); % top face center
    \coordinate (FrontC)  at ($ (A)!.5!(F) $); % front face center
    \coordinate (BackC)   at ($ (Dpt)!.5!(G) $); % back face center
    \coordinate (LeftC)   at ($ (A)!.5!(Hpt) $); % left face center
    \coordinate (RightC)  at ($ (B)!.5!(G) $); % right face center

    \foreach \q in {BottomC, LeftC, BackC} {
      \fill[fghblue!60, opacity=1.0] (\q) circle (5pt);
    }
    \foreach \q in {} {
      \fill[myred!60, opacity=1.0] (\q) circle (5pt);
    }

    % Transparent faces
    \fill[fghblue!20, opacity=0.25] (A) -- (B) -- (C) -- (Dpt) -- cycle;
    \fill[fghblue!20, opacity=0.25] (E) -- (F) -- (G) -- (Hpt) -- cycle;
    \fill[fghblue!20, opacity=0.25] (A) -- (B) -- (F) -- (E) -- cycle;
    \fill[fghblue!20, opacity=0.25] (Dpt) -- (C) -- (G) -- (Hpt) -- cycle;

    % Edges
    \draw[dashed, line width=2pt, fghblue] (A) -- (B);
    \draw[thick, line width=2pt, fghblue] (B) -- (C);
    \draw[thick, line width=2pt, fghblue] (C) -- (Dpt);
    \draw[dashed, line width=2pt, fghblue] (Dpt) -- (A);
    \draw[thick, line width=2pt, fghblue] (E) -- (F) -- (G) -- (Hpt) -- cycle;
    \draw[dashed, line width=2pt, fghblue] (A) -- (E);
    \draw[thick, line width=2pt, fghblue] (B) -- (F);
    \draw[thick, line width=2pt, fghblue] (C) -- (G);
    \draw[thick, line width=2pt, fghblue] (Dpt) -- (Hpt);

    \node[font=\footnotesize, anchor=south west] at (E) {$(1.1, 0.82, 1.2)$};
    \node[font=\footnotesize, anchor=north east] at (C) {$(0.9, 0.78, 0.6)$};

    % Corner points as little balls
    %\foreach \p in {A,B,C,Dpt,E,F,G,Hpt} {
    \foreach \p in {A,E,Hpt} {
      %\filldraw[fill=white, draw=black] (\p) circle (1.6pt);
      \fill[myred!60, opacity=1.0] (\p) circle (5pt);
    }
    \node[font=\footnotesize, anchor=south east] at (A) {$u_{25}$};
    \node[font=\footnotesize, anchor=south east] at (E) {$u_{27}$};
    %\node[font=\footnotesize, anchor=north west] at (G) {$u_{3}$};
    \node[font=\footnotesize, anchor=south west] at (Hpt) {$u_{21}$};
    
    \foreach \p in {B,C,Dpt,F,G} {
      \fill[fghblue!60, opacity=1.0] (\p) circle (5pt);
    }
    
    % Edge midpoints (centers of each edge)
    \coordinate (ABm) at ($ (A)!.5!(B) $);
    \coordinate (BCm) at ($ (B)!.5!(C) $);
    \coordinate (CDm) at ($ (C)!.5!(Dpt) $);
    \coordinate (DAm) at ($ (Dpt)!.5!(A) $);

    \coordinate (EFm) at ($ (E)!.5!(F) $);
    \coordinate (FGm) at ($ (F)!.5!(G) $);
    \coordinate (GHm) at ($ (G)!.5!(Hpt) $);
    \coordinate (HEm) at ($ (Hpt)!.5!(E) $);

    \coordinate (AEm) at ($ (A)!.5!(E) $);
    \coordinate (BFm) at ($ (B)!.5!(F) $);
    \coordinate (CGm) at ($ (C)!.5!(G) $);
    \coordinate (DPHm) at ($ (Dpt)!.5!(Hpt) $);

    % Edge centers as little balls
    %\foreach \q in {ABm, BCm, CDm, DAm, EFm, FGm, GHm, HEm, AEm, BFm, CGm, DPHm} {
    \foreach \q in {ABm, BCm, CDm, DAm, EFm, FGm, HEm, AEm, BFm, CGm, DPHm} {
      \fill[fghblue!60, opacity=1.0] (\q) circle (5pt);
    }
    \foreach \q in {GHm} {
      \fill[myred!60, opacity=1.0] (\q) circle (5pt);
    }
    %\node[font=\footnotesize, anchor=south west] at (FGm) {$u_{6}$}; % delete and add 12 and 23
    \node[font=\footnotesize, anchor=north west] at (GHm) {$u_{12}$}; % delete and add 12 and 23

    \foreach \q in {FrontC, RightC, TopC} {
      \fill[fghblue!60, opacity=1.0] (\q) circle (5pt);
    }
    \foreach \q in {} {
      \fill[myred!60, opacity=1.0] (\q) circle (5pt);
    }

    \coordinate (Center) at ($ (A)!.5!(G) $);
    \fill[fghblue!60, opacity= 1.0] (Center) circle (5pt);

    \coordinate (Nom) at ($ (BottomC)!.65!(TopC) $);
    \fill[black!60, opacity= 1.0] (Nom) circle (5pt);
    \node[font=\footnotesize, anchor=south] at (Nom) {$u_{nom}$};

    % Add a simple 3D coordinate system from the box origin (A)
    \def\L{2}
    \def\OrigX{4}
    \def\OrigY{6}
    \def\OrigZ{0}
    \coordinate (O) at (\OrigX,\OrigY,\OrigZ);L       % same as A
    \coordinate (Xaxis) at (\OrigX-\L,\OrigY,\OrigZ);  % end along x
    \coordinate (Yaxis) at (\OrigX,\OrigY-\L,\OrigZ);  % end along y
    \coordinate (Zaxis) at (\OrigX,\OrigY,\OrigZ+\L);  % end along z

    \draw[->, very thick, blue]  (O) -- (Xaxis) node[right] {$F_{12}$};
    \draw[->, very thick, red]   (O) -- (Yaxis) node[left]  {$w_{MF}$};
    \draw[->, very thick, purple] (O) -- (Zaxis) node[left] {$l$};
    
  \end{tikzpicture}
        \captionof{figure}{Reference discretization used for the MeOH-MF example. All scenarios that appear as worst-case in one of the scalarizations are marked in red. The nominal scenario is marked in black.}
        \label{fig:MeOHMFscenarios}
\end{figure}

\begin{table}[H]
    \centering
    \begin{tabular}{|c|c|c|c|c|}
        \hline
        \textbf{Scen.} & $\mathbf{F_{12}}$ & $\mathbf{w_{MF}}$ & $\mathbf{l}$ & \textbf{\replaced[id=R2]{worst-case}{WC} occurence} \\
        \hline \hline
%        3 & 0.9 & 0.78 & 1.2 & $w_{MeOH,bot}^{min}$ \\
%        \hline
        12 & 1.0 & 0.78 & 1.2 & $w_{MeOH,bot}^{min}$ \\
        \hline
        21 & 1.1 & 0.78 & 1.2 & $w_{MeOH,bot}^{min}$, $w_{MF,top}^{min}$ \\
        \hline
%        23 & 1.1 & 0.80 & 1.2 & $Q^{max}_c$, $Q^{max}_r$ \\
%        \hline
        25 & 1.1 & 0.82 & 0.6 & OPEX \\
        \hline
        27 & 1.1 & 0.82 & 1.2 & $Q^{max}_c$, $Q^{max}_r$, $F^{max}$ \\
        \hline
    \end{tabular}
    \captionof{table}{List of \replaced[id=R2]{worst-case}{WC} scenarios and the respective functions, in which each scenario leads to the worst-case. Some objectives and constraints have multiple worst-cases because these were determined across all Pareto points.}
    \label{tab:scenlist}
\end{table}

%\todo[inline, backgroundcolor=blue!20!white]{Comparison to nominal front and standard robust front}

\begin{figure}[H]
    \includegraphics[width=.75\textwidth]{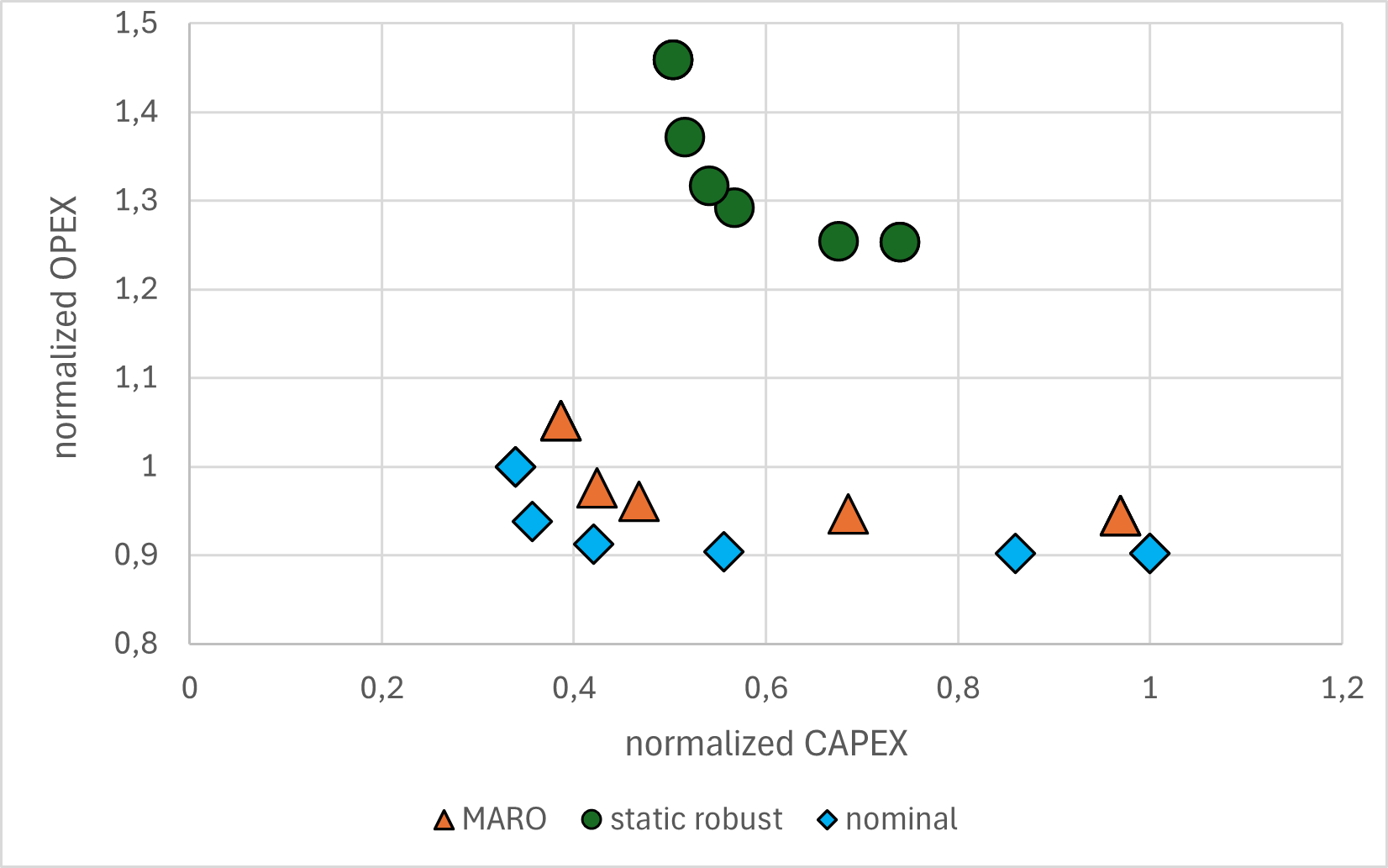}
    \caption{MARO Pareto front (orange triangle), nominal front (green circles) and non-adjustable robust front (blue diamonds) for the MeOH-MF separation problem. \added[id=R2]{CAPEX and OPEX have been normalized to the range of the nominal front.}}
    \label{fig:MeOHMFParetoPointsWStandardRob}
\end{figure}

Using the Sandwiching algorithm and the adaptive discretization scheme as described above, we obtain the MARO front of the MeOH-MF separation problem as depicted in Fig.~\ref{fig:MeOHMFParetoPointsWStandardRob}. As comparison, we additionally show both the nominal front and the non-adjustable robust front. For the non-adjustable robust front, we do not allow for adjustment of the wait-and-see decisions; see Ref.~\citenum{schwientek2025} for a more detailed discussion. As can be seen in Fig.~\ref{fig:MeOHMFParetoPointsWStandardRob}, the MARO Pareto front lies much closer to the nominal front than to the non-adjustable front. This illustrates the importance of taking the adjustability of the WSV into account for the MeOH-MF separation problem. Tab.~\ref{tab:MO_MRO_MARO_optimization_ranges} shows the ranges of the objectives, the HNV and the WSV for the nominal, the adjustable robust, and the non-adjustable robust MOs.

% Requires: \usepackage{amsmath}
\begin{table}[H]
    \centering
    \begin{tabular}{|c|c|c|c|c|}
        \hline
        & \textbf{MO} & \textbf{MARO} & \textbf{MRO} \\
        \hline \hline
        $ \text{(norm.) CAPEX}$ & 0.339 - 1.0 & 0.387 - 0.969 & 0.504 - 0.740 \\
        \hline
        $ \text{(norm.) OPEX}$ & 0.902 - 1.0 & 0.945 - 1.052 & 1.253 - 1.459\\
        \hline \hline
        $N \in [10, 150]$ & 12 - 150 & 13 - 150 & 19 - 61 \\
        \hline
        $N_f \in [3, 40]$ & 4 - 6 & 4 - 5 & 4 - 5 \\
        \hline
        $D \in [0.8, 2.0]$ & 0.955 - 1.006 &  1.071 - 1.133 & 1.257 - 1.367 \\
        \hline
        $A_r \in [50, 1000]$ & 150.80 - 167.13 & 189.46 - 210.85 & \replaced[id=R12]{251.01}{253.147} - 390.12 \\
        \hline
        $A_c \in [50, 1000]$ & 132.72 - 147.96 & 169.55 - 189.73 & 228.46 - 272.77 \\
        \hline \hline
        $R_V \in [0.5, 2.0]$ & 0.591 - 0.774 & 0.632 - 0.826 & 1.268 - 1.710 \\
        \hline
        $\dot{Q}_r \in [0.0625, 0.375]$ & 0.188 - 0.209 & 0.197 - 0.220 & 0.261 - 0.310 \\
        \hline
    \end{tabular}
    \caption{Ranges of (normalized) objectives and decision variables found in the MO, MARO, and MRO runs for the MeOH-MF separation. \added[id=R2]{CAPEX and OPEX have been normalized to the range of the nominal front.}}
    \label{tab:MO_MRO_MARO_optimization_ranges}
\end{table}

%\todo[inline, backgroundcolor=blue!20!white]{Robustification is necessary ("Gain of robustness")}

By hedging against all scenarios from the adaptive scheme, a degradation of obtainable OPEX and CAPEX is incurred for the MARO front compared to the nominal front. Conversely, this means that for the HNV obtained for each point of the nominal front, the WSV cannot be adjusted for every scenario in a way that retains both the achieved values for CAPEX and OPEX, and feasibility with respect to the problem constraints. We confirmed this for each of the \replaced[id=R12]{five}{eight} solutions of the nominal front approximation plotted in Fig.~\ref{fig:MeOHMFParetoPointsWStandardRob}. To this aim, for each scenario from the reference discretization, we re-optimized the WSV, while the HNV remained fixed. For any single one of the \replaced[id=R12]{five}{eight} nominal solutions, the re-optimization problems could find an adjustment of the WSV that was neither infeasible nor inferior for only \replaced[id=R12]{11-13}{13-15} out of the 28 scenarios (see Tab.~\ref{tab:scens_for_nom_front}).      
%This is illustrated in Fig.~\ref{}...
%\todo[inline]{Ich bin nicht ganz überzeugt, dass es sinnvoll ist, reoptimierte Punkte zu zeigen, wenn sie nicht zur Zulässigkeit optimiert wurden. Alternativ könnte man den Anteil der zur Zulässigkeit optimierbaren Szenarien einfach hier nennen.}
%\todo[inline, backgroundcolor=yellow!20!white]{Falls Bild gewünscht: hier einfügen. Vlt. hat MB schon ein entsprechendes Bild?}

%\color{red}{NOTE: The following table was computed anew.}
\begin{table}
    \centering
    \begin{tabular}{|l|ccccc|}
    \hline
        \textbf{nominal solution}  & \textbf{1} & \textbf{2} & \textbf{3} & \textbf{4} & \textbf{5} \\
        \hline
        \hline
        No. of scenarios with infeasible outcome & 7 & 13 & 13 & 15 & 7 \\
        \hline
        No. of scenarios with feasible but inferior outcome& 8 & 2 & 2 & 2 & 10 \\
    \hline
    \end{tabular}
    \caption{Accumulated outcomes when trying to adjust the WSV of the \replaced[id=R12]{five non-dominated nominal front solutions}{eight nominal front solutions} to each of the 28 reference scenarios. Solutions were numbered left to right according to Fig.~\ref{fig:MeOHMFParetoPointsWStandardRob}.}
    \label{tab:scens_for_nom_front}
\end{table}

%\todo[inline, backgroundcolor=blue!20!white]{Scenario fronts}

\begin{figure}[H]
    \includegraphics[width=.75\textwidth]{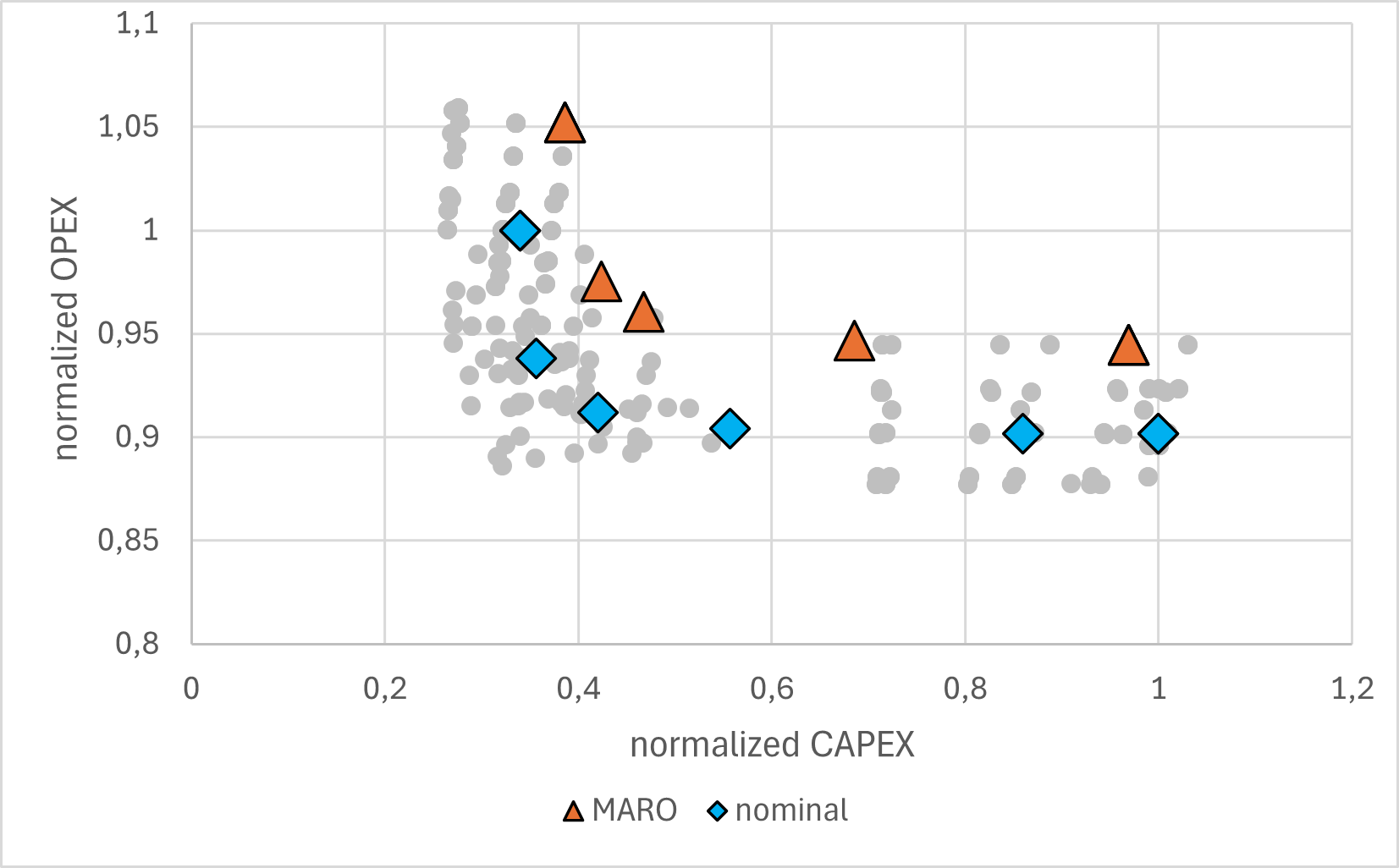}
    \caption{\replaced[id=R12]{Scenario front solutions}{Scenario fronts} for all scenarios from the reference discretization \replaced{(grey dots)}{(lines)} compared to the MARO front (orange triangles) \added[id=R12]{and the nominal front (blue diamonds)} for the MeOH-MF separation. \deleted[id=R12]{The apparent superiority of some MARO solutions is due to the pointwise approximation of the Pareto fronts. The MARO front was approximated more finely in the relevant area, whereas this was not the case for the scenario fronts, meaning that only the linear interpolation of the more widely spaced Pareto points there can be seen.}}
    \label{fig:scenario_fronts}
\end{figure}

Analogously to the nominal front, for each of the 28 scenarios from the reference discretization we can define the \emph{scenario front}: the set of Pareto optimal solutions where both HNV and WSV are optimized for the particular scenario. By definition, the MARO Pareto front is an upper bound to all scenario fronts. In Fig.~\ref{fig:scenario_fronts}, the approximation of the MARO front is plotted alongside the approximations of the 28 scenario fronts. For all front approximations, an adaptive Sandwiching algorithm was used.

%\todo[inline, backgroundcolor=blue!20!white]{Benefit of adaptive scheme (same front, less calculation time)}

\begin{figure}[H]
    \includegraphics[width=.75\textwidth]{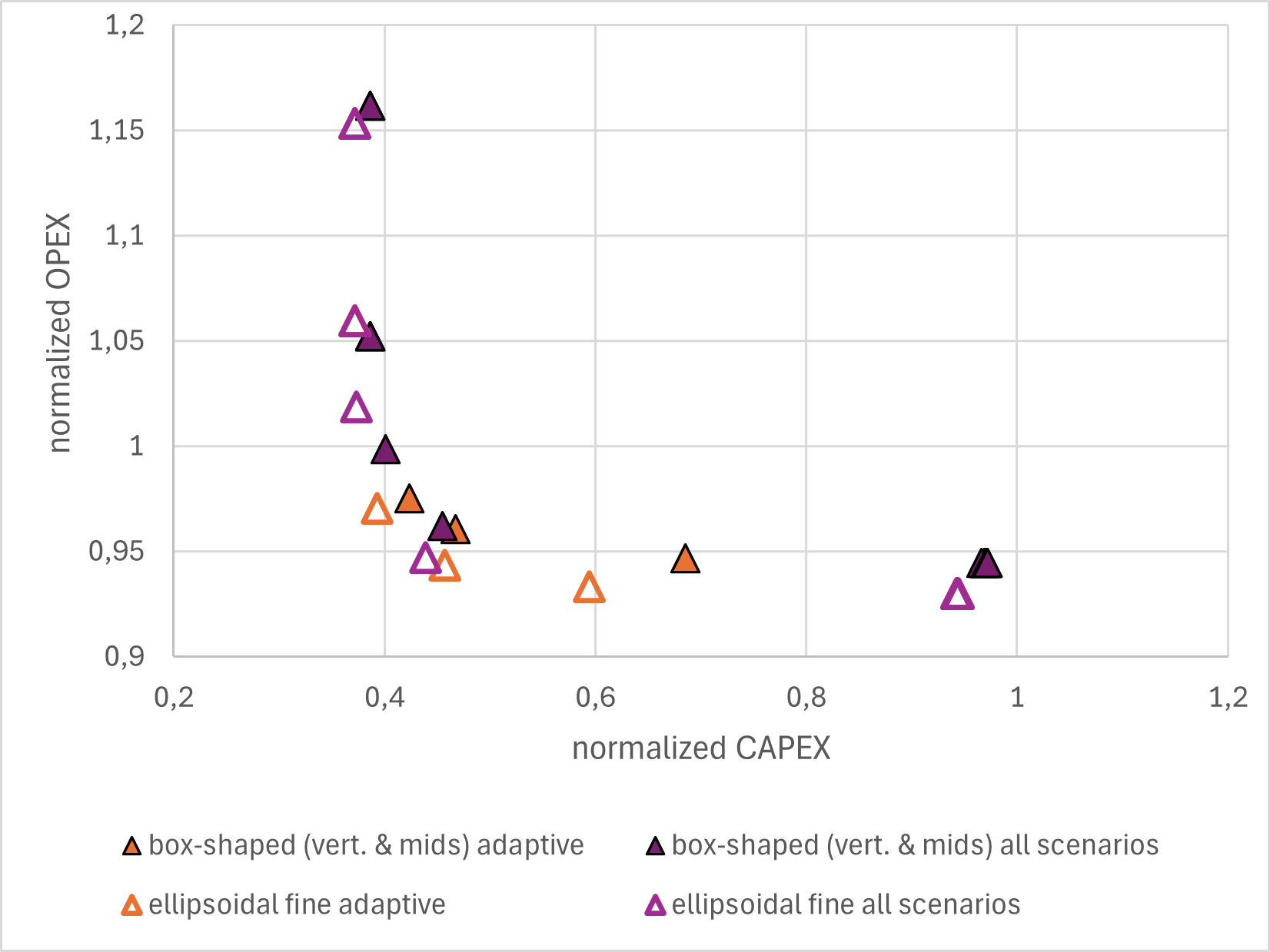}
    \caption{Both for the box-shaped with face mids (solid markers) and fine ellipsoidal (open markers) discretization, the MARO front for the adaptive scheme (orange) coincides with the front obtained when utilizing all 28 scenarios \replaced[id=R12]{(purple)}{(green circles)} for the MeOH-MF separation.}
    \label{fig:adaptive_vs_all}
\end{figure}

The number of optimization variables for the MARO problem ~(\ref{eq:discMARO}) depends on the number of scenarios, as for each scenario the WSV are replicated. In turn, the number of variables greatly impacts the calculation time. 
\replaced[id=R12]{Hence, by reducing the number of scenarios, the adaptive scenario selection algorithm is expected to significantly reduce the calculation time compared to using the full underlying reference discretization. Tab. \ref{tab:calculation_times} demonstrates this benefit for four underlying reference discretizations: box-shaped with only vertices, box-shaped with vertices and face mids, a coarse and a fine discretization on an ellipsoid. For each discretization, the MARO front calculation time is recorded using the adaptive choice versus using the full reference discretization. We see a reduction of calculation time between 79\% and 96\%. Fig.~\ref{fig:adaptive_vs_all} displays the calculated Pareto fronts for two of the discretization schemes, showing that the same fronts are obtained irrespective of the use of all scenarios or the adaptive scenario selection.}{We demonstrate the benefit of the adaptive scenario choice algorithm by solving the MARO problem (\ref{eq:discMARO}) using the unaltered reference discretization consisting of all 28 scenarios as a comparison. Using all scenarios, we obtain the same Pareto front (see Fig.~\ref{fig:adaptive_vs_all}), but the calculations took 767 minutes, which is a more than tenfold increase over the calculation time of 68 minutes when using the reduced scenario set from the adaptive scenario choice.}

%\color{red}{NOTE: The following table is new.}
\begin{table}
    \centering
    \begin{tabular}{|l|c|c|c|c|}
    \hline
        \textbf{}  & \makecell{\textbf{box-shaped} \\ \textbf{vertices only} \\ 9 scenarios} & \makecell{\textbf{box-shaped} \\ \textbf{vert. \& face mids} \\ 28 scenarios} & \makecell{\textbf{ellipsoidal} \\ \textbf{coarse} \\ 15 scenarios} & \makecell{\textbf{ellipsoidal} \\ \textbf{fine} \\ 28 scenarios}\\
        \hline
        \hline
        adaptive & 309 s & 1165 s & 597 s & 719 s \\
        \hline
        full discretization & 2494 s & 26767 s & 2824 s & 9486 s \\
%        \hline
%        comp. time saving & 88\% & 96\% & 79\% & 92\% \\
    \hline
    \end{tabular}
    \caption{Calculation times for four discretizations: two nested discretization for the uncertainty box and two nested discretization of an uncertainty ellipsoid with same extent as the box.}
    \label{tab:calculation_times}
\end{table}

\subsection{Price of Robustness}

%\todo[inline, backgroundcolor=blue!20!white]{Calculation of nominal scenario evaluations alongside approximation}

\begin{figure}[H]
    \includegraphics[width=.75\textwidth]{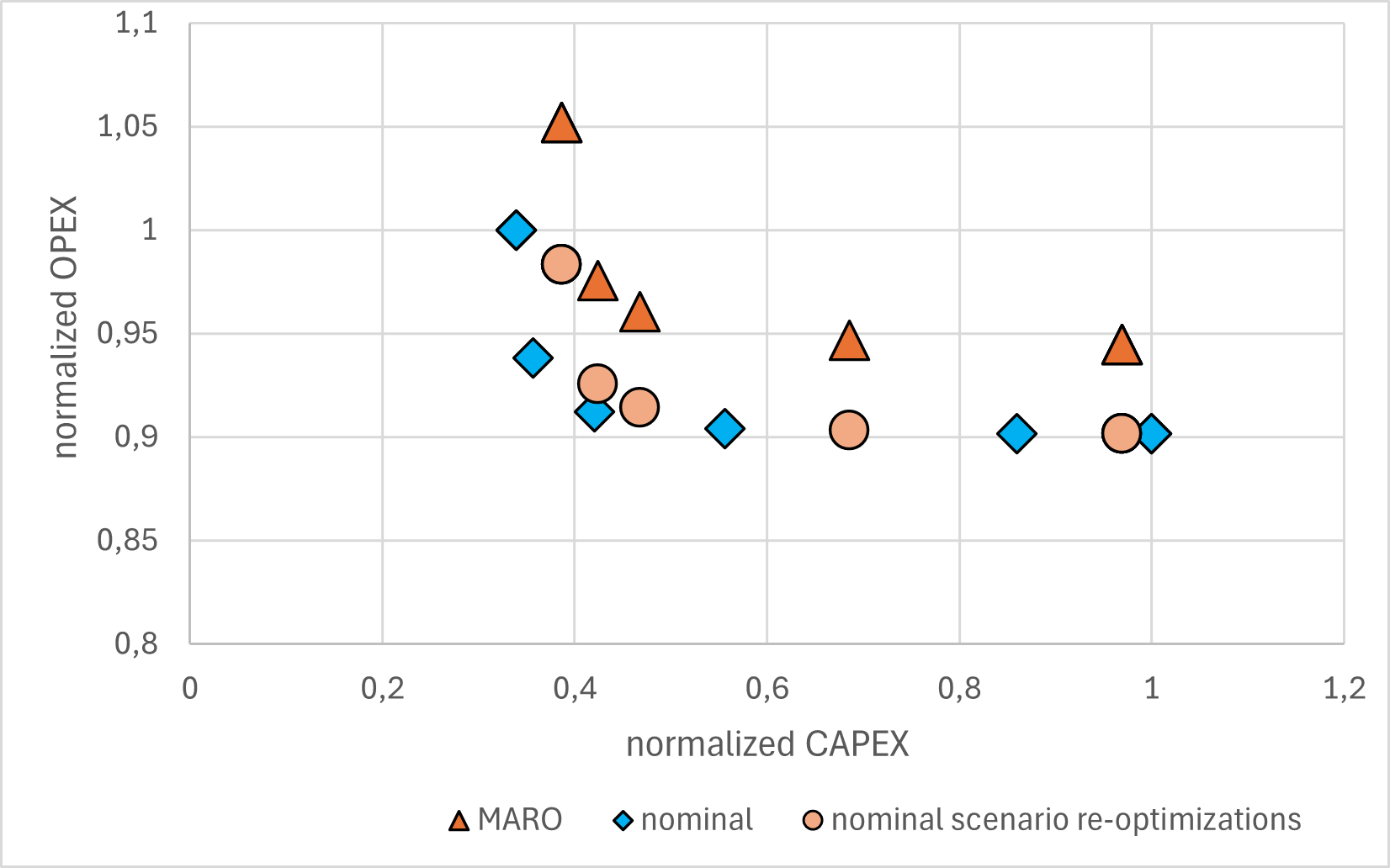}
    \caption{Nominal \replaced[id=R12]{(blue diamonds)}{(green circles)} and MARO (orange triangles) Pareto fronts as well as nominal scenario re-optimization points (orange circles) for the MeOH-MF separation.}
    \label{fig:MeOHMFParetoPoints}
\end{figure}

\begin{figure}[H]
    \includegraphics[width=.75\textwidth]{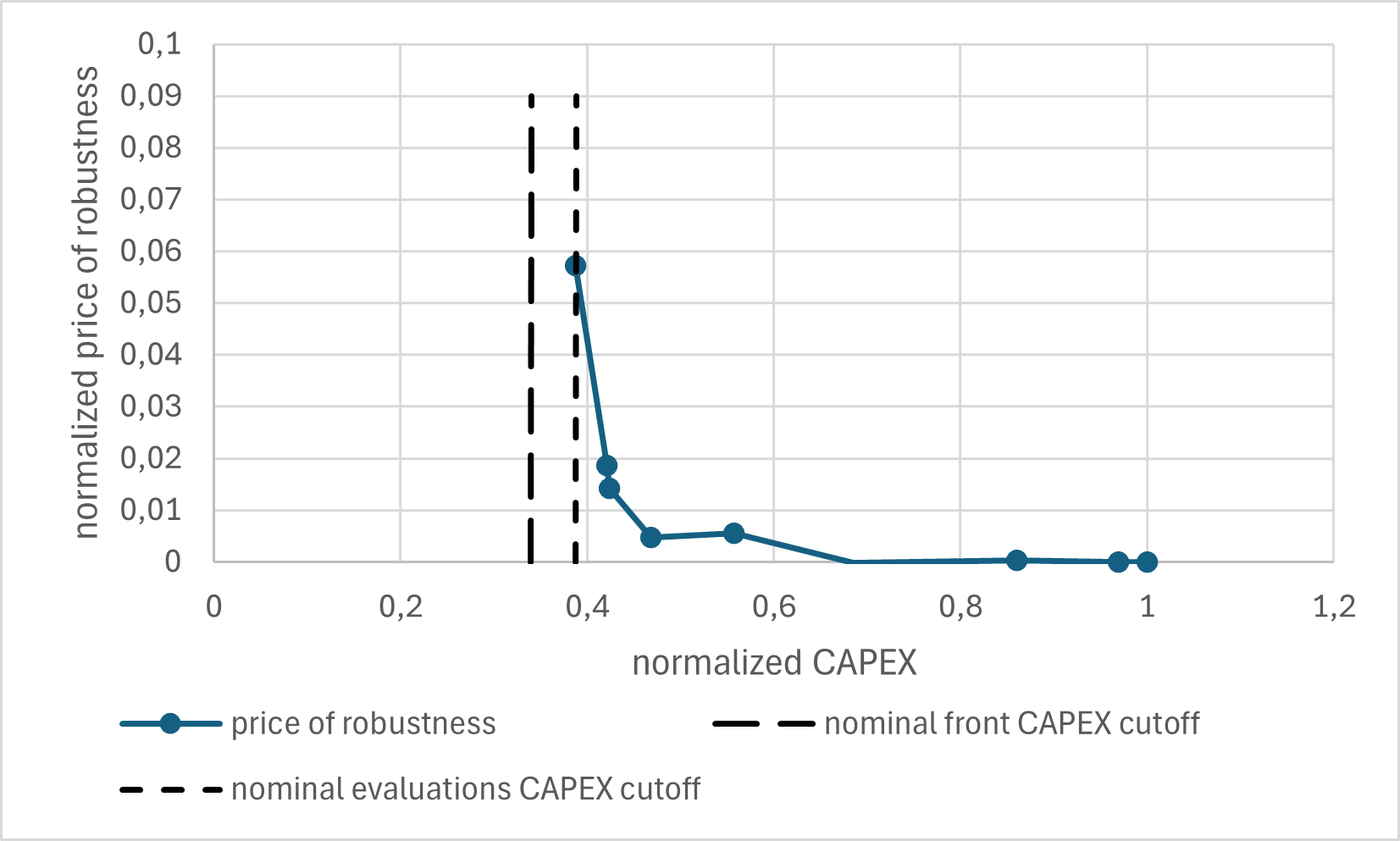}
    \caption{The price of robustness as a function of CAPEX for the MeOH-MF separation. For high CAPEX, the price of robustness is zero.}
    \label{fig:MeOHMF_price}
\end{figure}

To determine the price of robustness, the NSR problem \eqref{eq:NSRpt} is solved for each Pareto efficient HNV to obtain the optimal WSV for the nominal scenario. Note that for the MeOH-MF separation example, CAPEX does not depend on the WSV (but is fixed by the values of the HNV), and therefore problem \eqref{eq:NSRpt} has a single optimal solution that minimizes OPEX. Fig.~\ref{fig:MeOHMFParetoPoints} shows the MARO front, the nominal front, and the solutions to the nominal scenario re-optimization problems. Tab.~\ref{tab:optimization_ranges} shows the ranges of the objectives CAPEX and OPEX, the HNV and the WSV for all three solution sets. The price of robustness, as the difference in OPEX between the re-optimization points and the nominal front points, is displayed in Fig.~\ref{fig:MeOHMF_price}.
%\highlight[id=MW]{HNV and WSV}\todo{MW: vorher erklaeren/einfuehren, was hier die HNV und WSV sind} for the nominal and the MARO front as well as for the NSR point set. 

% Requires: \usepackage{amsmath}
\begin{table}[H]
    \centering
    \begin{tabular}{|c|c|c|c|c|}
        \hline
        & \textbf{MO} & \textbf{MARO} & \textbf{NSR} \\
        \hline \hline
        $ \text{(norm.) CAPEX}$ & 0.339 - 1.0 & 0.387 - 0.969 & 0.387 - 0.969 \\
        \hline
        $ \text{(norm.) OPEX}$ & 0.902 - 1.0 & 0.945 - 1.052 & 0.902 - 0.983 \\
        \hline \hline
        $N \in [10, 150]$ & 12 - 150 & 13 - 150 & 13 - 150 \\
        \hline
        $N_f \in [3, 40]$ & 4 - 6 & 4 - 5 & 4 - 5 \\
        \hline
        $D \in [0.8, 2.0]$ & 0.955 - 1.007 &  1.071 - 1.133 & 1.071 - 1.133 \\
        \hline
        $A_r \in [50, 1000]$ & 150.80 - 167.13 & 189.46 - 210.85 & 189.46 - 210.85 \\
        \hline
        $A_c \in [50, 1000]$ & 132.72 - 147.96 & 169.55 - 189.73 & 169.55 - 189.73 \\
        \hline \hline
        $R_V \in [0.5, 2.0]$ & 0.591 - 0.774 & 0.632 - 0.826 & 0.591 - 0.774 \\
        \hline
        $\dot{Q}_r \in [0.0625, 0.375]$ & 0.188 - 0.209 & 0.197 - 0.220 & 0.189 - \replaced[id=R12]{0.205}{0.209} \\
%        \hline \hline
%        Run time & ~5s per Pareto pt. & ~13min per Pareto pt. & - \\ % to be re-done, because of changes
        \hline
    \end{tabular}
    \caption{Ranges of (normalized) objectives and decision variables found in the optimization runs for the MeOH-MF separation.}
    \label{tab:optimization_ranges}
\end{table}

%\todo[inline, backgroundcolor=blue!20!white]{PoR over CAPEX discussion}

%The MARO Pareto front approximation, the nominal front approximation and the reoptimization points for all MARO Pareto optimal points are shown in Fig.~\ref{fig:MeOHMFParetoPoints}. The price of robustness, as the difference in OPEX between the reoptimization points and the nominal front points, is displayed in Figure \ref{fig:MeOHMF_price}. 

For high capital expenditure, nominal solutions and NSR points show equal normalized OPEX, i.e. the cost of robustness is zero. This shows that for a more expensive design of the distillation column, optimal operation in the nominal scenario can be obtained by adjusting the operating variables only. However, for lower capital expenditure, the operation variables cannot completely compensate for the suboptimal design in the nominal scenario. Also, the minimally possible capital expenditure for any feasible adjustable robust solution is markedly higher than for any feasible nominal solution.

%\todo[inline, backgroundcolor=green!20!white]{Als Linie; normalized CAPEX ab 0.2; markieren, bis wohin der Range in CAPEX für die nominale Front geht.}

%\todo[inline, backgroundcolor=blue!20!white]{PoR in navigation}

\begin{figure}[H]
	\includegraphics[width=.85\textwidth]{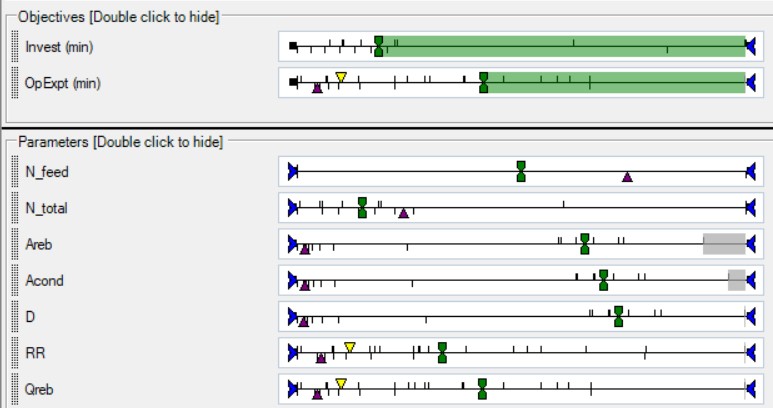}
	\caption{The interactive slider interface for the MeOH-MF separation problem. The currently selected MARO solution is represented by the green hourglasses. Yellow triangular markers depict the nominal scenario reoptimized WSV ("RR", "Qreb") and corresponding OPEX ("OpExpt"), while purple triangular markers show the nominal, non-robust solution.}
    \label{fig:MeOHMFsilders}
\end{figure}

Following the description in the section on real-time approximation and interactive display of the price of robustness, we have integrated the price of robustness into interactive navigation within the BASF software \textrm{CHEMASIM}. Fig.~\ref{fig:MeOHMFsilders} shows the slider interface for the MeOH-MF example. Sliders for the objectives (CAPEX as "Invest", and OPEX as "OpExpt") can be interactively dragged with the green hourglass handle. The yellow marker on the "OpExpt" slider signifies the value obtained by re-optimizing the WSV for the nominal scenario, while the purple marker shows the value on the nominal front. The difference between the two markers is the price of robustness. As the solution is changed by the engineer by dragging the hourglass handle, both markers update in real-time. 

The change in HNV and WSV can be observed in the sliders below the objectives ("Parameters"). The MARO variable values are shown as green hourglasses, and the nominal variable values as purple markers. For the WSV, namely the reflux ratio as "RR" and the heat duty as "Qreb", the reoptimized values are shown as yellow markers. 

\added[id=R2]{Within the decision making process, the display of the price of robustness allows the engineer to better adjust the level of robustness they are aiming for. If, for example, the price of robustness is too high, the engineer may decide to lower the required level of robustness by reducing the intervals of the uncertainties, or change the shape of the uncertainty set from box-shaped to ellipsoidal.}   

%\begin{figure}
%	\includegraphics[width=.45\textwidth]{gfx/priceofrobustness_standard.png}
%	\includegraphics[width=.45\textwidth]{gfx/priceofrobustness_adjustable.png}
%	\caption{Price of robustness for a standard robust Pareto set (lhs) and an adjustable robust Pareto set (rhs).}
%\end{figure}

\section{Conclusions and outlook}

In this work, adaptive scenario selection and computing the price of robustness were newly applied to model-based multi-criteria adjustable robust process design optimization.

Thinking of the realization of uncertainties as discrete scenarios turned out to be advantageous, both for computational reasons and from a user perspective. The worst-case scenarios for the determination of the worst-case Pareto solutions can be selected adaptively. This keeps the problem dimensions, and thus the computing times moderate without compromising the solution quality.

Hedging against uncertainties comes at a price, which has been introduced as the price of robustness. This concept may contribute to making discussions about safety margins more transparent, since the gain, i.e., increased robustness, can be traded against this price, which is a compromise in the performance indicators. 

So far, it has been assumed that the uncertainty set $U$ is given. From an application perspective, this is not necessarily the case. Rather, tolerance limits can be imposed on performance indicators $F$ in the MARO problem \refeq{eq:MARO}. This implies that an upper limit on the price of robustness is given. The question then is: Given this upper limit, which is the uncertainty set $U$ hedged against? Or, put as an optimization problem: Maximize the size, commonly the volume $\mathrm{Vol}$, of the uncertainty set $U$ that is hedged against while respecting an upper limit on the price of robustness.

This concept is known as inverse robustness \cite{berthold2024}. Within the multi-objective framework, maximizing $\mathrm{Vol}(U)$ is treated as an additional objective, while the size-determining parameters are additional decision variables. Once the Pareto-optimal solutions have been obtained, they can be interpolated as described above, thus smoothly interpolating the size of $U$.

For inverse-robust, multi-criteria process design optimization, it is also conceivable to select the scenarios adaptively as well as to quantify the robustness costs. For the former, it may even be necessary to refine the reference discretization due to the variability of the uncertainty set. This, along with the consideration of more complex uncertainty sets and discretizations, are topics for future research.

%The same framework enables inverse robustness: the range of scenarios covered by the operating conditions is added as an extra objective in the MO, yielding precise trade-offs between robustness and key performance indicators such as cost efficiency, yield, and product quality.

%Thinking in terms of what is and what is not affordable in terms of such a compromise leads to inverse robustness: If an effort is made to hedge against uncertainties, then as many scenarios shall be covered as possible. 

%In this work, rather simple scenario sets were dealt with. However, due to the highly non-linear relation between performance indicators and model parameters, it would be highly interesting to explore inverse robustness for more general sets. This remains a challenge for future work.

%In order to address inverse robustness, $U$ has to be parametrized in a way suitable for optimization. For the rectangular and elliptical sets in Fig. \ref{fig:discretizations}, the scaling factors in each dimension have been taken as optimization variables. In Fig.~\ref{fig:invRob}, results for a scaling the rectangular $U$-shape with equal scaling factors in each dimension are shown. 

%\begin{figure}
%	\includegraphics[width=.45\textwidth]{gfx/inverse_robustness_standard.png}
%	\includegraphics[width=.45\textwidth]{gfx/inverse_robustness_adjustable.png}
%	\caption{inverse robust Pareto sets in the standard robust (lhs) and adjustable robust approach (rhs). }
%   \label{fig:invRob}
%\end{figure}

\begin{acknowledgement}

The authors thank Kerstin Schneider and Elisabeth Halser for fruitful discussions. Part of this work was funded by the Deutsche Forschungsgemeinschaft (DFG, German Research Foundation) - GRK 2982, 516090167 "Mathematics of Interdisciplinary Multiobjective Optimization"

\end{acknowledgement}

%\nocite{*}
\bibliography{bibliography}

\end{document}